\documentclass[11pt,oneside]{amsart}
\usepackage{amsmath,amsthm,amssymb,mathtools}
\usepackage[dvipsnames,x11names,svgnames]{xcolor}
\usepackage[colorlinks=true,linkcolor=NavyBlue,citecolor=DarkGreen,urlcolor=blue]{hyperref}
\usepackage[margin=1in]{geometry}
\newcommand{\Fqx}{\mathbb{F}_q[x]}
\newcommand{\R}{\mathbb{R}}

\newcommand{\C}{\mathbb{C}}
\newcommand{\Fq}{\mathbb{F}_q}
\newcommand{\Hg}{\mathcal{H}_{2g+1}}
\renewcommand{\Re}{\textup{Re}}
\renewcommand{\Im}{\textup{Im}}
\newcommand{\legendre}[2]{\left(\frac{#1}{#2}\right)}
\newtheorem{theorem}{Theorem}
\newtheorem{lemma}[theorem]{Lemma}
\newtheorem*{thm A*}{Theorem A}
\newtheorem*{thm B*}{Theorem B}

\begin{document}

\title[]{Universality of $L$-functions in the Hyperelliptic Ensemble}

\author{Debanshu Ghosh}
\address{Department of Mathematics, University of Rochester, Rochester, NY 14627, USA}
\email{\url{dghosh5@ur.rochester.edu}}

\begin{abstract}
Let $\Fq$ be a finite field with $q$ elements, where $q$ is odd, and let $\Fqx$ be the polynomial ring over $\Fq$. We denote by $\Hg$ the set of monic, square-free polynomials $D \in \Fqx$ of degree $2g+1$. This is the hyperelliptic ensemble of the title. Associated with each $D$ is a non-trivial quadratic Dirichlet character $\chi_D$ and a corresponding Dirichlet $L$-function. We show that such $L$-functions are universal in the following sense. Let $U \subset \mathbb{C}$ be the open rectangular region with vertices $\sigma_1 \pm iA$ and $\sigma_2 \pm iA$, where $\tfrac12 < \sigma_1 < \sigma_2 < 1$ and $0 < A \le \frac{\pi}{3\log q}$, and let $K \subset U$ be a non-empty compact set. Let $F$ be a prescribed holomorphic, non-vanishing function on $U$ that is positive on $U \cap \mathbb{R}$. Then for a positive proportion of $D \in \Hg$, the function $L(s,\chi_D)$ approximates $F$ uniformly on $K$ arbitrarily closely.
\end{abstract}
\maketitle

\section{Introduction}\label{Sec Intro}

Let $s=\sigma+it$ and let $\zeta(s)$ denote the Riemann zeta function. In 1975, S.~M.~Voronin~\cite{vo 1} proved the following remarkable ``universality'' theorem for $\zeta(s)$.

\begin{thm A*}[Voronin]\label{thm A}
Let $D_r$ be the open disc of radius $r < \frac{1}{4}$ centered at $s = \frac{3}{4}$ in the complex plane. Suppose $F(s)$ is analytic on $D_r$ and continuous and non-vanishing on $\overline{D_r}$. Then for any $\varepsilon > 0$, there is a $\tau \in \mathbb{R}$ such that
\[
\max_{s \in \overline{D_r}} \left| \zeta(s + i\tau) - F(s) \right| < \varepsilon .
\]
\end{thm A*}

This theorem was soon extended in several directions by A.~Good~\cite{Good}, S.~M.~Gonek~\cite{gonekphd}, Bagchi~\cite{Bagchi 1,Bagchi 2}, Voronin himself~\cite{vo 2}--\cite{vo 5}, as well as others. For example, these authors proved results like Theorem~A for Dirichlet $L$-functions modulo $Q$, for the simultaneous approximation of finite sets of Dirichlet $L$-functions to various moduli, and for various Hurwitz zeta functions, while H.~Mishou and H.~Nagoshi~\cite{MN} established a universality result for quadratic Dirichlet $L$-functions $L(s,\chi_d)$. An analogue of Voronin's theorem has even been proved in the context of characteristic polynomials associated with unitary, orthogonal, and symplectic random matrices by Liza Jones~\cite{LJ}. In the years since, many papers studying different aspects of universality have appeared in the literature. Useful general surveys containing many of these results can be found in the paper by Matsumoto~\cite{Mmoto} and the two monographs by Laurinčikas~\cite{bookLaurinčikas} and Steuding~\cite{steuding2007value}.

For some years there has been considerable interest in extending number theoretical results from the number field setting to that of function fields. This was also recently done in the case of Voronin's Theorem~A by Andrade, Gonek, and Lee~\cite{AGL}. To describe their main theorem, we first require some notation concerning function fields and their $L$-functions. Excellent general references are Rosen~\cite{Rosen} and Thakur~\cite{Thakur}.

Let $\Fq$ be a finite field with $q$ elements, where $q$ is odd, and let $\Fqx$ denote the polynomial ring over $\Fq$. For any nonzero polynomial $f \in \Fqx$, the \emph{norm} of $f$ is defined by $|f| = q^{\deg f}$, and we set $|0| = 0$. When a monic polynomial $P$ is irreducible, we refer to it as a \emph{prime polynomial} or simply a \emph{prime}. For $\Re(s) > 1$, the zeta function of $\Fqx$ is defined by the series
\[
    \zeta_{\Fqx}(s) = \sum_{f \text{ monic}} \frac{1}{|f|^s}.
\]
Because $\Fqx$ contains exactly $q^j$ monic polynomials of degree $j$, we may group the terms by degree to express this sum as
\begin{equation}\label{eq_zeta_closed}
    \zeta_{\Fqx}(s) = \sum_{j=0}^{\infty} \frac{q^j}{q^{js}} = \frac{1}{1 - q^{1-s}}.
\end{equation}
Similarly, we can define the $L$-functions over $\Fqx$. Let $\chi$ be a non-principal Dirichlet character modulo $Q$. Then the associated $L$-function is initially defined for $\Re(s) > 1$ by the Euler product
\[
L(s,\chi) = \prod_{P \text{ prime}} \left( 1 - \chi(P)|P|^{-s} \right)^{-1}.
\]
This product expands directly into the Dirichlet series
\begin{equation}\label{eq_lsdef}
L(s,\chi) = \sum_{f \text{ monic}} \frac{\chi(f)}{|f|^s}, \qquad \Re(s) > 1.
\end{equation}

With these definitions in hand, the universality result of Andrade, Gonek, and Lee is the following.
\begin{thm B*}[Andrade--Gonek--Lee]\label{thm B}
Let $U$ be the open rectangle with vertices $\sigma_1+i\alpha$, $\sigma_2+i\alpha$, $\sigma_2+i\beta$, and $\sigma_1+i\beta$, where $\tfrac12<\sigma_1<\sigma_2<1$ and $0<\beta-\alpha\le\frac{2\pi}{3\log q}$. Let $C$ be a compact subset of $U$ with positive Lebesgue measure whose complement in $U$ is connected. Let $F(s)$ be a non-vanishing function that is analytic on the interior of $C$ and continuous on its boundary. Suppose also that $\mathcal{Q}$ is either the set of all primes or the set of polynomials that are products of a fixed finite set of primes. Then for any $\varepsilon>0$ and monic polynomial $Q\in\mathcal{Q}$ of sufficiently large degree, there exists a Dirichlet character $\chi \pmod Q$ such that
\[
\max_{s\in C}|L(s,\chi)-F(s)|<\varepsilon.
\]
Moreover, this approximation holds for a positive proportion of the $\Phi(Q)$ characters $\chi \pmod Q$, where $\Phi(Q)$ denotes the Euler totient function on $\Fqx$.
\end{thm B*}

Our goal in this paper is to prove an analogue of Theorem~B for quadratic $L$-functions in the hyperelliptic ensemble. In this way, our main result extends the work of Mishou and Nagoshi~\cite{MN} to the function-field setting. Before stating our result, we first introduce the definition of these $L$-functions.

Let $P \in \Fqx$ be a prime. Then by Proposition~1.10 of M.~Rosen~\cite{Rosen}, if $f \in \Fqx$ and $P \nmid f$, we know that the congruence $x^d \equiv f \pmod P$ is solvable if and only if
\[
f^{\frac{|P|-1}{d}} \equiv 1 \pmod P,
\]
where $d$ is a divisor of $q-1$. The quadratic residue symbol arises when we consider $d=2$ and is denoted by $\Bigl( \dfrac{f}{P} \Bigr) \in \{1,-1\}$, satisfying
\[
\Bigl( \dfrac{f}{P} \Bigr) \equiv f^{\frac{|P|-1}{2}} \pmod P
\]
for $f$ coprime to $P$. Thus, if $P$ is a prime and $f \in \Fqx$, the quadratic residue symbol $\Bigl( \dfrac{f}{P} \Bigr)$ can be written as
\[
\Bigl( \dfrac{f}{P} \Bigr) =
\begin{cases}
  0, & \text{if } P \mid f,\\[6pt]
  1, & \text{if } P \nmid f \text{ and } f \text{ is a square mod } P,\\[6pt]
 -1, & \text{if } P \nmid f \text{ and } f \text{ is a non-square mod } P.
\end{cases}
\]
If $Q = P_1^{e_1} P_2^{e_2} \cdots P_k^{e_k}$ is the prime factorization of a monic polynomial $Q \in \Fqx$, we extend this to the \emph{Jacobi symbol} defined as
\[
\Bigl( \dfrac{f}{Q} \Bigr) = \prod_{i=1}^k \Bigl( \dfrac{f}{P_i} \Bigr)^{e_i}.
\]
If $c \in \Fq^\ast \, (= \Fq \setminus \{0\})$, then
\[
\Bigl( \dfrac{c}{Q} \Bigr) = c^{\frac{(q-1)}{2}\deg Q},
\]
For coprime monic polynomials $A,B \in \Fqx$, the law of quadratic reciprocity, proved by Artin~\cite{Artin}, asserts that
\[
\Bigl( \dfrac{A}{B} \Bigr)
= \Bigl( \dfrac{B}{A} \Bigr)
(-1)^{\frac{(q-1)}{2}\deg A \deg B},
\]
and if $A,B$ share a common factor, both sides vanish.

Given a monic, square-free polynomial $D \in \Fqx$, we now define the quadratic character $\chi_D$ for any monic $f \in \Fqx$ by
\[
\chi_D(f) = \Bigl( \dfrac{D}{f} \Bigr).
\]
Let $\Hg = \{ D \in \Fqx : D \text{ monic, square-free, } \deg D = 2g + 1 \}$. Then each polynomial $D \in \Hg$ determines a quadratic Dirichlet character $\chi_D$, and each such character $\chi_D$ gives rise to the Dirichlet $L$-series
\[
  L(s,\chi_D) = \sum_{f \text{ monic}} \frac{\chi_D(f)}{|f|^s},
  \qquad \Re(s) > 1,
\]
which in product form is
\[
L(s,\chi_D) = \prod_{P \text{ prime}} (1 - \chi_D(P)|P|^{-s})^{-1},
\qquad \Re(s) > 1.
\]
It is often more natural to work with the variable $u = q^{-s}$. In terms of $u$, these two equivalent representations become
\[
\mathcal{L}(u,\chi_D) = \sum_{f \text{ monic}} \chi_D(f) u^{\deg f}, \qquad |u| < \frac{1}{q},
\]
and
\[
\mathcal{L}(u,\chi_D) = \prod_{P \text{ prime}} \left( 1 - \chi_D(P) u^{\deg P} \right)^{-1}, \qquad |u| < \frac{1}{q}.
\]
Note that $\mathcal{L}(u,\chi_D)$ is a polynomial of degree $2g$ (see Rosen~\cite{Rosen}, Proposition~4.3). Moreover, it is well known that it satisfies a Riemann hypothesis (see Weil~\cite{Aweil}), namely, all of its zeros lie on the circle $|u| = q^{-\frac12}$. Hence, letting $e(x) = e^{2\pi i x}$, one may factor $\mathcal{L}(u,\chi_D)$ as
\[
\mathcal{L}(u,\chi_D) = \prod_{j=1}^{2g} (1 - \alpha_j u),
\]
where $\alpha_j = q^{\frac12} e(-\theta_j)$, $j = 1, \dots, 2g$, are the reciprocals of the roots $u_j = q^{-\frac12} e(\theta_j)$ of $\mathcal{L}(u,\chi_D)$.

We may now state our universality theorem for hyperelliptic ensembles.

\begin{theorem}\label{thm 1}
Let $U \subset \C$ denote the open rectangular region with vertices $\sigma_1 \pm iA$ and $\sigma_2 \pm iA$, where
\[
\frac12 < \sigma_1 < \sigma_2 < 1 \quad \text{and} \quad 0 < A \le \frac{\pi}{3\log q}.
\]
Let $K \subset U$ be a non-empty compact set, and let $F$ be a non-vanishing holomorphic function on $U$ that is positive on $U \cap \R$. Then for any $\varepsilon > 0$, the inequality
\[
\max_{s \in K} |L(s, \chi_D) - F(s)| < \varepsilon
\]
holds for a positive proportion of $D \in \Hg$ as $g \to \infty$, where $\chi_D$ is the quadratic character defined above and $L(s,\chi_D)$ is its associated Dirichlet $L$-function. More precisely, we have
\[
\liminf_{g \to \infty} \frac{1}{\#\Hg} \# \Big\{ D \in \Hg : \max_{s \in K} |L(s, \chi_D) - F(s)| < \varepsilon \Big\} > 0.
\]
\end{theorem}

\section{The Hybrid Approximation Formula}

To establish Theorem~\ref{thm 1}, we make use of a truncated product formula for $L(s, \chi_D)$. Andrade, Gonek, and Keating established the following lemma. For the proof, see~\cite{AGK}.

\begin{lemma}\label{lem_hybrid_product}
Let $g \ge 1$ and $D \in \Hg$. For any integer $k \ge 0$, we define the truncated Euler factor $\mathcal{P}_k(u, \chi_D)$ by
\begin{equation}\label{eq_def_pk}
    \mathcal{P}_k(u, \chi_D) = \exp\Bigg( \sum_{n=1}^k \, \sum_{\substack{f \mathrm{\ monic} \\ \deg f = n}} \frac{\Lambda(f) \chi_D(f)}{n} u^n \Bigg),
\end{equation}
where $\Lambda(f)$ is the von Mangoldt function on $\Fqx$, given by
\[
    \Lambda(f) = \begin{cases}
        \deg P, & \text{if } f = P^n \text{ for some prime } P \text{ and } n \ge 1, \\
        0, & \text{otherwise}.
    \end{cases}
\]
We define the corresponding analytic tail factor $\mathcal{Z}_k(u, \chi_D)$ by
\[
    \mathcal{Z}_k(u, \chi_D) = \exp\Bigg( -\sum_{j=1}^{2g} \, \sum_{n > k} \frac{(\alpha_j u)^n}{n} \Bigg).
\]
Then for all $|u| < q^{-\frac12}$, the function $\mathcal{L}(u, \chi_D)$ admits the exact representation
\begin{equation}\label{eq_trunc_inu}
    \mathcal{L}(u, \chi_D) = \mathcal{P}_k(u, \chi_D) \mathcal{Z}_k(u, \chi_D).
\end{equation}
\end{lemma}

A useful approximation for this hybrid formula in terms of $s$ is given by the following lemma.

\begin{lemma}\label{lem_trunc_approx}
Let $g \ge 1$, $D \in \Hg$, $\sigma_1 > \frac12$, and let $k \ge 1$ be an integer such that $\frac{g}{k} q^{\left(\frac12 - \sigma_1\right)k} \le 1$. Then uniformly for all $s \in \C$ with $\Re(s) \ge \sigma_1$, we have
\begin{equation}\label{eq_hybrid_exp}
    L(s, \chi_D) = P_k(s, \chi_D) \left( 1 + O \left( \frac{g}{k} q^{\left(\frac12 - \sigma_1\right)k} \right) \right),
\end{equation}
where
\begin{equation}\label{eq_def_pk_s}
    P_k(s, \chi_D) = \exp\Bigg( \sum_{\substack{P \mathrm{\ prime } \\ \deg P \le k}} \sum_{n=1}^{\lfloor k / \deg P \rfloor} \frac{\chi_D(P^n)}{n |P|^{ns}} \Bigg).
\end{equation}
The implied constant in the $O$-term depends at most on $q$ and $\sigma_1$.
\end{lemma}

\begin{proof}
Let $u = q^{-s}$. For $\Re(s) \ge \sigma_1 > \frac12$, we have $|u| \le q^{-\sigma_1} < q^{-\frac12}$. Grouping the sum in~\eqref{eq_def_pk} over prime powers $f = P^n$, we see that the condition $\deg f \le k$ becomes $n \deg P \le k$, which gives
\begin{align*}
    \sum_{r=1}^k \sum_{\substack{f \text{ monic} \\ \deg f = r}} \frac{\Lambda(f) \chi_D(f) q^{-rs}}{r}
    &= \sum_{\substack{P \text{ prime} \\ \deg P \le k}} \sum_{n=1}^{\lfloor k / \deg P \rfloor} \frac{(\deg P)\, \chi_D(P^n)}{n \deg P} q^{-n s \deg P} \\
    &= \sum_{\substack{P \text{ prime} \\ \deg P \le k}} \sum_{n=1}^{\lfloor k / \deg P \rfloor} \frac{\chi_D(P^n)}{n |P|^{ns}},
\end{align*}
which establishes~\eqref{eq_def_pk_s}.
Next, we bound the tail factor $\log \mathcal{Z}_k(q^{-s}, \chi_D)$. Since $|\alpha_j| = q^{\frac12}$ and $|u| \le q^{-\sigma_1}$, the triangle inequality gives
\begin{align*}
    |\log \mathcal{Z}_k(q^{-s}, \chi_D)|
    &\le \sum_{j=1}^{2g} \sum_{n=k+1}^{\infty} \frac{|\alpha_j|^n |u|^n}{n} \\
    &\le 2g \sum_{n=k+1}^{\infty} \frac{q^{n\left(\frac12 - \sigma_1\right)}}{n} \\
    &\le \frac{2g}{k+1} q^{k\left(\frac12 - \sigma_1\right)} \frac{q^{\frac12 - \sigma_1}}{1 - q^{\frac12 - \sigma_1}} \\
    &\ll_{q,\sigma_1} \frac{g}{k} q^{\left(\frac12 - \sigma_1\right)k}.
\end{align*}
Since $\frac{g}{k} q^{\left(\frac12 - \sigma_1\right)k} \le 1$, applying the estimate $e^w = 1 + O(w)$ on bounded discs and substituting the above bound into~\eqref{eq_trunc_inu} proves~\eqref{eq_hybrid_exp}.
\end{proof}

\section{Approximation and Density Lemmas}\label{Sec ApproxLemmas}

In this section, we introduce the analytic and approximation lemmas required to prove the Fundamental Lemma (Lemma~\ref{lem_fundamental}) in Section~\ref{Sec ProofLemma}. Specifically, we establish two discretization results for replacing real weights in $[-1, 1]$ by signs in $\{1, -1\}$, prove an $L^2$-density theorem for Dirichlet polynomials with real weights, and record a real-coefficient variant of Runge's theorem. We emphasize to the reader that the development in this section and the next closely follows Sections~4 and~5 of~\cite{AGL}, with suitable modifications to fit our setting.

Throughout this section and Section~\ref{Sec ProofLemma}, let $U$ denote the open rectangular region with vertices $\sigma_1 \pm iA$ and $\sigma_2 \pm iA$, where $\tfrac{1}{2} < \sigma_1 < \sigma_2 < 1$ and $0 < A \le \frac{\pi}{3\log q}$. Further, let $K \subset U$ be a non-empty compact set.

\begin{lemma}\label{lem_linear_algebra}
Let $M$ and $N$ be positive integers with $M \le N$. Let $a_{m,n}$ and $b_m$, for $1 \le m \le M$ and $1 \le n \le N$, be real numbers. Suppose the system of equations
\begin{equation}\label{eq_linequ}
\sum_{n=1}^{N} a_{m,n} x_n = b_m \qquad (1 \le m \le M)
\end{equation}
has a solution $\textbf{x} = (x_1,\ldots,x_N) \in D_N$, where $D_N = \{(x_1,\ldots,x_N) : x_n \in [-1,1], \, 1 \le n \le N\}$. Then the system also has a solution $\textbf{x}' \in D_N$ for which $x_n' = \pm 1$ for at least $N-M$ integers $n$.
\end{lemma}

\begin{proof}
We proceed by induction on $N$ for a fixed $M$. The base case $N = M$ is trivial. For $N > M$, the solutions to~\eqref{eq_linequ} define a linear manifold of dimension at least $N-M$ over $\mathbb{R}$. Since this manifold intersects the convex body $D_N$, it must intersect its boundary. Thus, there exists a solution with at least one coordinate $x_k = \pm 1$.

Fixing such a coordinate reduces~\eqref{eq_linequ} to a system of $M$ equations in $N-1$ variables. After reordering the variables, we have
\[
\sum_{n=1}^{N-1} a_{m,n} x_n + a_{m,N} x_N = b_m \qquad (1 \le m \le M).
\]
By the inductive hypothesis, this reduced system has a solution with at least $(N-1)-M$ coordinates of absolute value $1$. Including the previously fixed coordinate, we obtain at least $N-M$ such coordinates in total. This proves the lemma.
\end{proof}

\begin{lemma}\label{lem_Lambdap_to_ap}
Let $K$ have its usual meaning defined above, and let $1 \le \ell < m$. Suppose that for each prime $P$ with $\ell < \deg P \le m$ we have a real number $\lambda_P \in [-1,1]$. Then there are numbers $a_P \in \{1,-1\}$ such that
\[
\max_{s \in K}
\Bigg|
\sum_{\substack{P \mathrm{\ prime }\\\ell < \deg P \le m}} \lambda_P |P|^{-s}
-
\sum_{\substack{P \mathrm{\ prime }\\\ell < \deg P \le m}} a_P |P|^{-s}
\Bigg|
\ll q^{-\frac{\ell}{2}},
\]
where the implied constant depends at most on $K$, $q$, and $\sigma_1$.
\end{lemma}

\begin{proof}
Let $B = 1 + \max_{s\in K} |s|$, so that $B$ depends only on $K$, $B > 1$, and $|s| < B$ for all $s \in K$. Set $R = \lfloor 2Be^2 k \log q \rfloor + 1$ with $k \ge 1$, and suppose $k < k' \le 2k$. Then for $s \in K$, we have
\begin{equation}\label{eq_S1+S2}
\begin{split}
 \sum_{\substack{P \text{ prime}\\ k < \deg P \le k'}} \lambda_P |P|^{-s} &= \sum_{\substack{P \text{ prime}\\ k < \deg P \le k'}} \lambda_P \left(\sum_{r=0}^{\infty} \frac{(-s \deg P \log q)^r}{r!}\right) \\
 &= \sum_{\substack{P \text{ prime}\\ k < \deg P \le k'}} \lambda_P \left\{ \sum_{r=0}^{R} + \sum_{r > R} \right\} \frac{(-s \deg P \log q)^r}{r!}.
\end{split}
\end{equation}

Now let $\pi_q(j)$ denote the number of primes $P$ of degree $j$. The prime number theorem for function fields asserts that
\begin{equation}\label{eq_PNT}
\pi_q(j) = \frac{q^j}{j} + O\Big(\frac{q^{\frac{j}{2}}}{j}\Big),
\end{equation}
and in particular $\pi_q(j) \le \frac{q^j}{j}$ for all $j \ge 1$. Thus, if $N$ denotes the number of primes $P$ with $k < \deg P \le k'$, then since $q \ge 3$ and $k' \le 2k$,
\begin{equation}\label{eq_no_of_primes}
N = \sum_{k < j \le k'} \pi_q(j) \ll \sum_{j=k}^{k'} \frac{q^j}{j} \ll q^{k'} \ll q^{2k}.
\end{equation}

To estimate the sum over $r > R$ in~\eqref{eq_S1+S2}, first observe that $e^r > \frac{r^r}{r!}$ for $r \ge 1$, so that $r! > (\frac{r}{e})^r$ for these $r$. Also, recall that $\lambda_P \in [-1,1]$, $|s| < B$ for $s \in K$, and $R = \lfloor 2Be^2 k \log q \rfloor + 1 > 2Be^2 k \log q$ with $B > 1$. Then by~\eqref{eq_no_of_primes}, the sum over $r > R$ in~\eqref{eq_S1+S2} satisfies
\begin{equation}\label{eq_S2}
\begin{split}
&\le \sum_{\substack{P \text{ prime}\\k < \deg P \le k'}} \sum_{r > R} \frac{(|s| \deg P \log q)^r}{r!} \\
&\ll q^{2k} \sum_{r > R} \left( \frac{Be \cdot 2k \log q}{r} \right)^r \le q^{2k} \sum_{r > R} \left( \frac{2Be k \log q}{R} \right)^r \\
&< q^{2k} \sum_{r > R} e^{-r} \ll q^{2k} e^{-R} \le q^{2k} e^{-2Be^2 k \log q} \\
&\le q^{2k - 2k e^2} < q^{-12k},
\end{split}
\end{equation}
with absolute implied constants (since $B > 1$ and $2 - 2e^2 < -12$). Hence we may write~\eqref{eq_S1+S2} as
\begin{equation}\label{eq_main+error_1}
\begin{split}
\sum_{\substack{P \text{ prime}\\k < \deg P \le k'}} \lambda_P |P|^{-s}
&= \sum_{\substack{P \text{ prime}\\k < \deg P \le k'}} \lambda_P \left( \sum_{r=0}^{R} \frac{(-s \deg P \log q)^r}{r!} \right) + O(q^{-12k}) \\
&= \sum_{r=0}^{R} s^r \Bigg(\sum_{\substack{P \text{ prime}\\k < \deg P \le k'}} \lambda_P \frac{(-\deg P \log q)^r}{r!} \Bigg) + O(q^{-12k}).
\end{split}
\end{equation}

To treat the double sum on the right, we consider separately the cases where $N$, the number of primes $P$ in the inner sum, satisfies $N \ge R+1$ or $N < R+1$. First consider the case where $N \ge R+1$. Set
\[
b_r = \sum_{\substack{P \text{ prime}\\k < \deg P \le k'}} \lambda_P \frac{(-\deg P \log q)^r}{r!} \qquad (0 \le r \le R).
\]
Then the system of equations
\[
\sum_{\substack{P \text{ prime}\\k < \deg P \le k'}} x_P \frac{(-\deg P \log q)^r}{r!} = b_r \qquad (0 \le r \le R)
\]
has a solution $x_P = \lambda_P$. Hence, by Lemma~\ref{lem_linear_algebra}, there exist real numbers $\lambda_P' \in [-1,1]$ such that $\lambda_P' \in (-1, 1)$ for at most $R+1$ of the primes $P$ and
\[
b_r = \sum_{\substack{P \text{ prime}\\k < \deg P \le k'}} \lambda_P' \frac{(-\deg P \log q)^r}{r!} \qquad (0 \le r \le R).
\]
Hence, by~\eqref{eq_main+error_1},
\[
\sum_{\substack{P \text{ prime}\\k < \deg P \le k'}} \lambda_P |P|^{-s} = \sum_{\substack{P \text{ prime}\\k < \deg P \le k'}} \lambda_P' \left(\sum_{r=0}^{R} \frac{(-s \deg P \log q)^r}{r!}\right) + O(q^{-12k}).
\]
As in~\eqref{eq_S2}, we can extend the sum over $r$ to infinity with an error of size $O(q^{-12k})$ (again with an absolute implied constant). In this way we see that
\[
\sum_{\substack{P \text{ prime}\\k < \deg P \le k'}} \lambda_P |P|^{-s} = \sum_{\substack{P \text{ prime}\\k < \deg P \le k'}} \lambda_P' |P|^{-s} + O(q^{-12k}).
\]
As we saw above, all but at most $R+1$ of the values $\lambda_P'$ equal $+1$ or $-1$. Thus, there are integers $a_P \in \{1,-1\}$ for which
\begin{equation}\label{eq_final_sum}
\begin{split}
\sum_{\substack{P \text{ prime}\\k < \deg P \le k'}} \lambda_P |P|^{-s} &= \sum_{\substack{P \text{ prime}\\k < \deg P \le k'}} a_P |P|^{-s} + O(q^{-12k}) + O(R q^{-k\sigma_1}) \\
&= \sum_{\substack{P \text{ prime}\\k < \deg P \le k'}} a_P |P|^{-s} + O(q^{-\frac{k}{2}}).
\end{split}
\end{equation}
Here we have used the fact that for $s \in K$, we have $\Re(s) > \sigma_1 > \frac12$, thus $R q^{-k\sigma_1} \ll_{K, q, \sigma_1} q^{-\frac{k}{2}}$. For the remaining case $N < R+1$, for any choice of $a_P \in \{1,-1\}$, the sums on both sides of~\eqref{eq_final_sum} are trivially bounded by
\[
\ll R q^{-k\sigma_1} \ll_{K, q, \sigma_1} q^{-\frac{k}{2}}.
\]

Thus, whether $N \ge R+1$ or $N < R+1$, we have for all $s \in K$ that
\begin{equation}\label{eq_diff_of_sums}
\bigg| \sum_{\substack{P \text{ prime}\\k < \deg P \le k'}} \lambda_P |P|^{-s} - \sum_{\substack{P \text{ prime}\\k < \deg P \le k'}} a_P |P|^{-s} \bigg| \ll q^{-\frac{k}{2}},
\end{equation}
where the implied constant depends at most on $K$, $q$, and $\sigma_1$. Now let $J = \lfloor \log_2(m/\ell) \rfloor \ge 0$, so that $2^J \ell \le m < 2^{J+1}\ell$. Applying~\eqref{eq_diff_of_sums} to each dyadic interval $(2^j\ell, 2^{j+1}\ell]$ for $0 \le j \le J - 1$ and once more to the remaining interval $(2^J \ell, m]$ (if $2^J \ell < m$), and summing the resulting errors, we find that
\[
\Bigg| \sum_{\substack{P \text{ prime}\\\ell < \deg P \le m}} \lambda_P |P|^{-s} - \sum_{\substack{P \text{ prime}\\\ell < \deg P \le m}} a_P |P|^{-s} \Bigg|
\ll \sum_{j=0}^{J} q^{-2^{j-1}\ell}
\ll q^{-\frac{\ell}{2}}
\]
uniformly for $s \in K$. This completes the proof.
\end{proof}

Before proceeding to the next result, we introduce the underlying Hilbert space. Let $U \subset \mathbb{C}$ be the open rectangular region bounded by the vertices $\sigma_1 \pm iA$ and $\sigma_2 \pm iA$, where
\[
\frac12 < \sigma_1 < \sigma_2 < 1
\qquad\text{and}\qquad
0 < A \le \frac{\pi}{3\log q}.
\]
We denote by $L^2(U)$ the space of complex-valued measurable functions on $U$ that are square-integrable with respect to the Lebesgue measure. On $L^2(U)$, we consider the real inner product given by
\[
\langle f,g\rangle = \Re \iint_U f(s)\overline{g(s)}\,d\sigma\,dt.
\]
The induced norm is defined by
\[
\|f\| = \sqrt{\langle f,f\rangle} = \left( \iint_U |f(s)|^2\,d\sigma\,dt \right)^{\frac12}.
\]
Equipped with this inner product and norm, $L^2(U)$ is a real Hilbert space. Let $P^2(U) \subset L^2(U)$ denote the closed subspace spanned by polynomials with real coefficients (equivalently, polynomials that are real-valued on the interval $U \cap \mathbb{R}$).

For $0 < \ell < m$, define
\[
T_{\ell,m} = \Bigg\{ \sum_{\substack{P \text{ prime}\\\ell < \deg P \le m}} \lambda_P |P|^{-s} : \lambda_P \in [-1,1] \Bigg\}.
\]
The sets $T_{\ell,m}$ are convex and nested, satisfying $T_{\ell,m_1} \subseteq T_{\ell,m_2}$ whenever $\ell < m_1 < m_2$. Thus, their union
\[
T_\ell = \bigcup_{m > \ell} T_{\ell,m}
\]
is also convex. Let $\overline{T_\ell}$ denote the closure of $T_\ell$ in $L^2(U)$. Since every function in $T_\ell$ is entire and real-valued on $U \cap \mathbb{R}$, its Taylor series around any point in $U \cap \mathbb{R}$ has real coefficients and converges uniformly on $\overline{U}$, so $T_\ell \subseteq P^2(U)$. Since $P^2(U)$ is closed, it follows that $\overline{T_\ell} \subseteq P^2(U)$. It remains to show the reverse inclusion, namely $P^2(U) \subseteq \overline{T_\ell}$.

\begin{lemma}\label{lem_functional_analysis}
For every $\ell \ge 1$, we have $P^2(U) \subseteq \overline{T_\ell}$.
\end{lemma}

\begin{proof}
Since $P^2(U)$ is a closed subspace of the Hilbert space $L^2(U)$, it is isomorphic to its dual. The set $T_\ell$ is convex, and consequently so is its closure $\overline{T_\ell}$. By a consequence of the Hahn--Banach theorem (see Rudin~\cite{WRudinFA}, Theorem~3.4), an element $g \in P^2(U)$ belongs to $\overline{T_\ell}$ if and only if
\[
\langle g, f \rangle \in \langle \overline{T_\ell}, f \rangle \quad \text{for all } f \in P^2(U).
\]
If $f \equiv 0$, this condition is trivially satisfied. For any $f \not\equiv 0$, it therefore suffices to show that $\langle \overline{T_\ell}, f \rangle = \mathbb{R}$.

To this end, we assume that $f \not\equiv 0$ and consider
\[
\langle T_{\ell,m}, f \rangle = \Bigg\{ \sum_{\substack{P \text{ prime}\\\ell < \deg P \le m}} \lambda_P \Re \iint_U \overline{f(s)} |P|^{-s} \, d\sigma \, dt : \lambda_P \in [-1, 1] \Bigg\},
\]
which is a closed symmetric interval in $\mathbb{R}$ centered at the origin with radius
\begin{equation}\label{eq_radi}
r_f(\ell, m) = \sum_{\substack{P \text{ prime}\\\ell < \deg P \le m}} \left| \Re \iint_U \overline{f(s)} |P|^{-s} \, d\sigma \, dt \right|.
\end{equation}
To demonstrate that these nested intervals cover $\mathbb{R}$, we seek a sequence $m_j \to \infty$ such that $r_f(\ell, m_j) \to \infty$ as $j \to \infty$. For this purpose, define for $x \in \mathbb{R}$
\[
G(x) = \iint_{U} \overline{f(s)} e^{-sx} \, d\sigma \, dt.
\]
Because $U$ is symmetric across the real axis, $x \in \mathbb{R}$, and $f \in P^2(U)$ satisfies $f(\bar{s}) = \overline{f(s)}$, the integrand $h(s) = \overline{f(s)} e^{-sx}$ satisfies $h(\bar{s}) = \overline{h(s)}$. Splitting the vertical integral over $[-A, 0]$ and $[0, A]$ and substituting $t \mapsto -t$ on $[-A, 0]$ gives
\begin{align*}
G(x) &= \int_{\sigma_1}^{\sigma_2} \int_{0}^{A} \left( \overline{f(\sigma - it)} e^{-(\sigma - it)x} + \overline{f(\sigma + it)} e^{-(\sigma + it)x} \right) dt \, d\sigma \\
&= \int_{\sigma_1}^{\sigma_2} \int_{0}^{A} \left( \overline{\overline{f(s)} e^{-sx}} + \overline{f(s)} e^{-sx} \right) dt \, d\sigma \\
&= 2 \Re \int_{\sigma_1}^{\sigma_2} \int_{0}^{A} \overline{f(s)} e^{-sx} \, dt \, d\sigma.
\end{align*}
Thus $G(x)$ is real-valued for all $x \in \mathbb{R}$, and
\[
\Re \iint_{U} \overline{f(s)} e^{-sx} \, d\sigma \, dt = \iint_{U} \overline{f(s)} e^{-sx} \, d\sigma \, dt.
\]

Now, for $\delta_0 \in (0, 1 - \sigma_2)$, define $H(x)$ for $x \in \mathbb{R}$ by
\[
H(x) = \Re \iint_{U} \overline{f(s)} e^{(\sigma_2 + \delta_0 - s)x} \, d\sigma \, dt = e^{(\sigma_2 + \delta_0)x} \iint_{U} \overline{f(s)} e^{-sx} \, d\sigma \, dt.
\]
Setting $x = \deg P \log q$, we see from~\eqref{eq_radi} that
\[
r_f(\ell, m_j) = \sum_{\substack{P \text{ prime}\\\ell < \deg P \le m_j}} q^{-(\sigma_2 + \delta_0)\deg P} \bigl| H(\deg P \log q) \bigr|.
\]
Since the number of primes $P$ of degree $d$ satisfies $\pi_q(d) \gg \frac{q^d}{d}$, grouping the sum by degree leads to the lower bound
\[
r_f(\ell, m_j) \gg \sum_{\ell < d \le m_j} \frac{q^{(1 - \sigma_2 - \delta_0)d}}{d} \bigl| H(d \log q) \bigr|.
\]
Now choose $\delta_1$ such that $0 < \delta_1 < 1 - \sigma_2 - \delta_0$. If there exists a sequence of integers $m_j \to \infty$ satisfying
\[
\bigl| H(m_j \log q) \bigr| > q^{-\delta_1 m_j},
\]
then
\begin{equation}\label{eq_r_to_infty}
r_f(\ell, m_j) \gg \frac{q^{m_j(1 - \sigma_2 - \delta_0 - \delta_1)}}{m_j} \xrightarrow{j \to \infty} \infty.
\end{equation}

Suppose, to the contrary, that this claim is false. Then there exists an integer $d_0 \ge 1$ such that, for all $d \ge d_0$,
\begin{equation}\label{eq_H_d_less_than}
\left|H(d\log q)\right| \le q^{-\delta_1 d}.
\end{equation}
Define the power series
\[
\widetilde{H}(z) = \sum_{d=0}^{\infty} H(d\log q) z^d.
\]
Since~\eqref{eq_H_d_less_than} holds for all sufficiently large $d$, the root test shows that $\widetilde{H}(z)$ converges absolutely and hence defines an analytic function in the open disc $D_1 = \{ z \in \C : |z| < q^{\delta_1} \}$. On the other hand, if $|z| < q^{-(\sigma_2-\sigma_1+\delta_0)}$, then
\begin{align*}
\widetilde{H}(z)
&= \sum_{d=0}^{\infty} \left( \iint_U \overline{f(s)} q^{(\sigma_2+\delta_0-s)d} \, d\sigma \, dt \right) z^d \\
&= \iint_U \overline{f(s)} \sum_{d=0}^{\infty} \left( z q^{\sigma_2+\delta_0-s} \right)^d \, d\sigma \, dt \\
&= \iint_U \frac{\overline{f(s)}}{1 - z q^{\sigma_2+\delta_0-s}} \, d\sigma \, dt.
\end{align*}

This integral representation defines an analytic function on the open set
\[
\Omega = \mathbb{C} \setminus \left\{ q^{-(\sigma_2 + \delta_0 - s)} : s \in \overline{U} \right\}.
\]
The excluded set is the compact annular sector $\{ z \in \C : q^{-(\sigma_2 - \sigma_1 + \delta_0)} \le |z| \le q^{-\delta_0}, \, |\arg z| \le A \log q \}$, which lies strictly inside $D_1$ (since $q^{-\delta_0} < 1 < q^{\delta_1}$) and has angular width $2A\log q \le \frac{2\pi}{3} < 2\pi$. Consequently, both $\Omega$ and the intersection $D_1 \cap \Omega$ are connected. Since the power series and the integral representation agree on the open disc $|z| < q^{-(\sigma_2 - \sigma_1 + \delta_0)} \subset D_1 \cap \Omega$, the identity theorem implies that they agree throughout $D_1 \cap \Omega$, and in particular on the annulus $q^{-\delta_0} < |z| < q^{\delta_1}$. For $|z| > q^{-\delta_0}$ and $s \in U$, we have $|z q^{\sigma_2+\delta_0-s}| > 1$. Using the identity $(1 - w)^{-1} = -\sum_{d=1}^{\infty} w^{-d}$ valid for $|w| > 1$, we obtain
\[
\widetilde{H}(z) = -\sum_{d=1}^{\infty} z^{-d} \iint_{U} \overline{f(s)} q^{-d(\sigma_2 + \delta_0 - s)} \, d\sigma \, dt.
\]
Thus, in the annulus $q^{-\delta_0} < |z| < q^{\delta_1}$, we obtain two Laurent series expansions for $\widetilde{H}(z)$, namely
\[
\sum_{d=0}^{\infty} H(d \log q) z^d = - \sum_{d=1}^{\infty} z^{-d} \iint_U \overline{f(s)} q^{-d(\sigma_2 + \delta_0 - s)} \, d\sigma \, dt.
\]
By the uniqueness of Laurent series expansions on an annulus, the corresponding coefficients on both sides must agree, meaning that each coefficient must vanish identically. Consequently,
\begin{equation}\label{eq_q_ds}
\iint_{U} \overline{f(s)} q^{ds} \, d\sigma \, dt = 0 \quad \text{for all } d \in \mathbb{Z}.
\end{equation}

Our next goal is to show that
\begin{equation}\label{eq_s_n}
\iint_{U} \overline{f(s)} s^n \, d\sigma \, dt = 0 \quad \text{for all integers } n \ge 0.
\end{equation}

For any $s = \sigma + it \in \overline{U}$, write $q^{s-1} = r e^{i\theta}$. Then $r \le q^{\sigma_2 - 1} < 1$ and $|\theta| = |t|\log q \le A\log q \le \frac{\pi}{3}$. The condition $|\theta| \le \frac{\pi}{3}$ implies $2\cos\theta \ge 1 > r$, which gives
\[
|1 - q^{s-1}|^2 = 1 - r(2\cos\theta - r) < 1.
\]
By the continuity of $|1 - q^{s-1}|$ on the compact set $\overline{U}$, it follows that
\begin{equation}\label{eq_M_bound}
\max_{s \in \overline{U}} |1 - q^{s-1}| = M < 1.
\end{equation}
For every integer $k \ge 0$, we have
\begin{equation}\label{eq_binomial_step}
\iint_U \overline{f(s)}(1 - q^{s-1})^k \, d\sigma \, dt = \sum_{j=0}^k \binom{k}{j}(-q^{-1})^j \left( \iint_U \overline{f(s)}q^{js} \, d\sigma \, dt \right) = 0,
\end{equation}
where the last equality follows from~\eqref{eq_q_ds}. Now since $\log(1 - w) = -\sum_{k=1}^\infty \frac{w^k}{k}$ for $|w| < 1$, we see that if $w = 1 - q^{s-1}$, then
\[
-\iint_U \overline{f(s)} \bigg(\sum_{k=1}^\infty \frac{1}{k} (1 - q^{s-1})^k \bigg)\, d\sigma \, dt = \iint_U \overline{f(s)}\log(q^{s-1}) \, d\sigma \, dt.
\]
By~\eqref{eq_M_bound}, we may interchange the order of summation and integration on the left-hand side and use~\eqref{eq_binomial_step} to obtain
\[
\iint_U \overline{f(s)}\log(q^{s-1}) \, d\sigma \, dt = 0.
\]
Hence,
\[
-\log q \iint_U \overline{f(s)} \, d\sigma \, dt + \log q \iint_U \overline{f(s)} s \, d\sigma \, dt = 0.
\]
By~\eqref{eq_q_ds} with $d = 0$, the first integral on the left-hand side vanishes. Because $\log q \ne 0$, the second integral must vanish as well. We therefore find that
\[
\iint_U \overline{f(s)} s^d \, d\sigma \, dt = 0
\]
for $d = 0, 1$.

To extend this to larger values of $d$, we consider the expansion
\[
\bigl(\log(q^{s-1})\bigr)^n = (-1)^n \sum_{k_1, \dots, k_n = 1}^\infty \frac{1}{k_1 \cdots k_n}(1 - q^{s-1})^{k_1 + \dots + k_n}
\]
for $n \ge 1$. We multiply both sides of this by $\overline{f(s)}$ and integrate over $U$. By~\eqref{eq_M_bound}, we may integrate the expression on the right term by term. Doing so and applying~\eqref{eq_binomial_step}, we find that
\begin{align*}
\iint_U \overline{f(s)}\bigl(\log(q^{s-1})\bigr)^n \, d\sigma \, dt &= (-1)^n \sum_{k_1, \dots, k_n = 1}^\infty \frac{1}{k_1 \cdots k_n} \iint_U \overline{f(s)}(1 - q^{s-1})^{k_1 + \dots + k_n} \, d\sigma \, dt \\
&= 0.
\end{align*}
Consequently,
\[
\iint_U \overline{f(s)}(-\log q + s\log q)^n \, d\sigma \, dt = 0.
\]
Expanding via the binomial theorem, we see that
\[
\sum_{d=0}^n \binom{n}{d} (-\log q)^{n-d}(\log q)^d \iint_U \overline{f(s)} s^d \, d\sigma \, dt = 0.
\]
By induction on $n$, we now deduce that
\[
\iint_U \overline{f(s)} s^n \, d\sigma \, dt = 0
\]
for all integers $n \ge 0$. This establishes~\eqref{eq_s_n}.

Finally, we show that $f \equiv 0$. Because $f \in P^2(U)$, there exists a sequence of polynomials $\{f_n(s)\}$ with real coefficients such that $\|f - f_n\| \to 0$ as $n \to \infty$. By~\eqref{eq_s_n},
\[
\iint_U \overline{f(s)} f_n(s) \, d\sigma \, dt = 0 \quad \text{for all } n \ge 1.
\]
By the Cauchy--Schwarz inequality, we have
\begin{align*}
\iint_U |f(s)|^2 \, d\sigma \, dt &= \iint_U \overline{f(s)} \bigl(f(s) - f_n(s)\bigr) \, d\sigma \, dt \\
&\le \|f\| \|f - f_n\| \xrightarrow{n \to \infty} 0.
\end{align*}
Thus $\|f\| = 0$, which implies that $f \equiv 0$, since $f$ is analytic inside $U$. This contradicts our earlier assumption that $f$ is not identically zero inside $U$. This shows that~\eqref{eq_H_d_less_than} is false, which in turn means that~\eqref{eq_r_to_infty} is true and completes the proof.
\end{proof}

We shall require below a variant of Runge's theorem producing polynomials with real coefficients. Before establishing that, we first recall the polynomial case of Runge's theorem (see Rudin~\cite[Theorem~13.9]{WRudin}).

\begin{lemma}[Runge's theorem]\label{lem_Runge}
Let $\Omega \subset \mathbb{C}$ be an open set whose complement in the Riemann sphere $\widehat{\mathbb{C}} \setminus \Omega$ is connected, and let $\tilde{K} \subset \Omega$ be a compact set. If $F(s)$ is holomorphic on $\Omega$ and $\varepsilon > 0$, then there exists a polynomial $p(s) \in \mathbb{C}[s]$ such that
\[
\max_{s \in \tilde{K}} |F(s) - p(s)| < \varepsilon.
\]
\end{lemma}

The variant of Runge's theorem we require is the following.
\begin{lemma}\label{lem_real_runge}
Let $U$ and $K$ be defined as at the beginning of Section~\ref{Sec ApproxLemmas}. Let $F$ be a holomorphic function on $U$ that is real-valued on $U \cap \mathbb{R}$. Then for every $\varepsilon > 0$, there exists a polynomial $p(s) \in \mathbb{R}[s]$ such that
\[
\max_{s \in K} |F(s) - p(s)| < \varepsilon.
\]
\end{lemma}

\begin{proof}
Since $K \subset U$ is a non-empty compact subset of the open rectangular region $U$, we may choose real numbers $\sigma_1', \sigma_2', A'$ satisfying
\[
\sigma_1 < \sigma_1' < \min_{s \in K} \Re(s) \le \max_{s \in K} \Re(s) < \sigma_2' < \sigma_2
\qquad\text{and}\qquad
\max_{s \in K} |\Im(s)| < A' < A.
\]
Let $U'$ be the open rectangular region with vertices $\sigma_1' \pm iA'$ and $\sigma_2' \pm iA'$. Then $\overline{U'}$ is a compact subset of $U$ that is symmetric about $\mathbb{R}$ and contains $K$.

Because $U$ is open, connected, and symmetric about $\mathbb{R}$ with $U \cap \mathbb{R} = (\sigma_1, \sigma_2) \neq \emptyset$, and $F$ is holomorphic on $U$ and real-valued on $U \cap \mathbb{R}$, the Schwarz reflection principle implies that
\[
F(\overline{s}) = \overline{F(s)} \quad \text{for all } s \in U.
\]
Since $U$ is an open rectangle, its complement $\widehat{\mathbb{C}} \setminus U$ in the Riemann sphere is connected. Applying Runge's theorem (Lemma~\ref{lem_Runge}) to the compact set $\overline{U'} \subset U$, we obtain a polynomial $p_0(s) = \sum_{j=0}^n c_j s^j \in \mathbb{C}[s]$ such that
\[
\max_{s \in \overline{U'}} |F(s) - p_0(s)| < \varepsilon.
\]
Define $p(s) \in \mathbb{R}[s]$ by averaging $p_0(s)$ with its conjugate polynomial $p_0^*(s) = \overline{p_0(\overline{s})}$. That is, define
\[
p(s) = \frac{p_0(s) + p_0^*(s)}{2} = \sum_{j=0}^n \Re(c_j) s^j.
\]
For every $s \in \overline{U'}$, the symmetry of $\overline{U'}$ implies $\overline{s} \in \overline{U'}$. Using the reflection identity $F(s) = \overline{F(\overline{s})}$, we obtain
\begin{align*}
|F(s) - p(s)| &\le \frac12|F(s) - p_0(s)| + \frac12|F(s) - \overline{p_0(\overline{s})}| \\
&= \frac12|F(s) - p_0(s)| + \frac12|\overline{F(\overline{s})} - \overline{p_0(\overline{s})}| \\
&< \varepsilon
\end{align*}
for all $s \in \overline{U'}$. Since $K \subset \overline{U'}$, this completes the proof.
\end{proof}

\begin{lemma}\label{lem_Titchmarsh}
Let $W$ be a bounded region in $\mathbb{C}$. Let $K$ be a compact set in $\mathbb{C}$ such that $K \subset W$. Suppose that $f(s)$ is analytic on $W$ and that
\[
\iint_W |f(s)|^2 \, d\sigma \, dt \le L,
\]
where $L > 0$. Then
\[
\max_{s \in K} |f(s)| \le a(W,K) L^{\frac12}.
\]
Here $a(W,K)$ is a positive constant depending only on $W$ and $K$.
\end{lemma}
A proof of essentially the same lemma can be found in Titchmarsh~\cite{Titchmarsh} (p.~303).

\section{Proof of the Fundamental Lemma}\label{Sec ProofLemma}

We now combine the approximation and density results from Section~\ref{Sec ApproxLemmas} to establish the Fundamental Lemma (Lemma~\ref{lem_fundamental}), which plays a central role in the proof of Theorem~\ref{thm 1} and adapts Lemma~15 of~\cite{AGL} to our setting. As a first step, we prove a preliminary version (Lemma~\ref{lem_baby_fundamental}) that approximates a holomorphic function $F$ on $K$ by a finite sum over primes whose degrees lie in an interval $(\ell, m]$.

\begin{lemma}\label{lem_baby_fundamental}
Let $U$ and $K$ be defined as at the beginning of Section~\ref{Sec ApproxLemmas}. Let $F$ be a holomorphic function on $U$ that is real-valued on $U \cap \R$. Then, for each $\ell \ge 1$, there exists a constant $m_0 = m_0(q, U, K, \ell, F) > \ell$ such that for every $m \ge m_0$, there exists a choice of numbers $a_P \in \{1,-1\}$ indexed over all primes $P \in \Fqx$ with $\ell < \deg P \le m$ satisfying
\begin{equation}\label{eq_fundamental_1_est}
\max_{s \in K} \Bigg| F(s) - \sum_{\substack{P \mathrm{\ prime }\\\ell < \deg P \le m}} \frac{a_P}{|P|^s} \Bigg| \le C q^{-\frac{\ell}{2}}.
\end{equation}
The constant $C$ depends only on $q$, $K$, and $U$.
\end{lemma}

\begin{proof}
Let $\ell \ge 1$ be given. By Lemma~\ref{lem_real_runge}, there exists a polynomial $p(s) \in \mathbb{R}[s]$ such that
\begin{equation}\label{eq_runge_approx}
\max_{s \in K} |F(s) - p(s)| < q^{-\frac{\ell}{2}}.
\end{equation}
Because $p(s)$ is a polynomial with real coefficients, we have $p(s) \in P^2(U)$. By Lemma~\ref{lem_functional_analysis}, $P^2(U) \subseteq \overline{T_\ell}$, which implies that $p(s) \in \overline{T_\ell}$. Consequently, there exist a real number $m_0 > \ell$, depending on $q$, $U$, $K$, $\ell$, and $p(s)$ (and hence on $F$), and real coefficients $\lambda_P \in [-1, 1]$ for all primes $P$ with $\ell < \deg P \le m_0$, such that
\[
\iint_{U} \Bigg| \sum_{\substack{P \text{ prime}\\\ell < \deg P \le m_0}} \lambda_P |P|^{-s} - p(s) \Bigg|^2 d\sigma \, dt \le q^{-\ell}.
\]
For any $m \ge m_0$, setting $\lambda_P = 0$ for all primes $P$ with $m_0 < \deg P \le m$ gives real coefficients $|\lambda_P| \le 1$ for all primes $P$ with $\ell < \deg P \le m$ satisfying
\[
\iint_{U} \Bigg| \sum_{\substack{P \text{ prime}\\\ell < \deg P \le m}} \lambda_P |P|^{-s} - p(s) \Bigg|^2 d\sigma \, dt \le q^{-\ell}.
\]
Because both $p(s)$ and the Dirichlet polynomial $\sum_{\ell < \deg P \le m} \lambda_P |P|^{-s}$ are analytic on $U$, Lemma~\ref{lem_Titchmarsh} gives
\begin{equation}\label{eq_lambda_p_to_ps}
\max_{s \in K} \Bigg| \sum_{\substack{P \text{ prime}\\\ell < \deg P \le m}} \lambda_P |P|^{-s} - p(s) \Bigg| \le a(U, K) q^{-\frac{\ell}{2}},
\end{equation}
where $a(U, K)$ depends only on $U$ and $K$.

Next, by Lemma~\ref{lem_Lambdap_to_ap}, the coefficients $\lambda_P \in [-1, 1]$ can be replaced by numbers $a_P \in \{1,-1\}$ so that
\begin{equation}\label{eq_sign_approx}
\max_{s \in K} \Bigg| \sum_{\substack{P \text{ prime}\\\ell < \deg P \le m}} \lambda_P |P|^{-s} - \sum_{\substack{P \text{ prime}\\\ell < \deg P \le m}} a_P |P|^{-s} \Bigg| \ll q^{-\frac{\ell}{2}},
\end{equation}
where the implied constant depends at most on $K$, $q$, and $\sigma_1$ (and hence on $q$, $K$, and $U$). Combining inequalities~\eqref{eq_runge_approx}, \eqref{eq_lambda_p_to_ps}, and~\eqref{eq_sign_approx} by the triangle inequality, we deduce~\eqref{eq_fundamental_1_est} with a constant $C$ depending only on $q$, $K$, and $U$.
\end{proof}

\begin{lemma}[The Fundamental Lemma]\label{lem_fundamental}
Let $U$, $K$, and $F$ be defined as in Lemma~\ref{lem_baby_fundamental}. Then, for any given integer $\ell_0 \ge 1$ and any real number $\varepsilon > 0$, there exist an integer $m = m(q, U, K, \varepsilon, F, \ell_0) > \ell_0$ and a choice of numbers $a_P \in \{1,-1\}$ for each prime $P$ with $\deg P \le m$, such that
\[
\max_{s \in K} \Bigg| F(s) - \sum_{\substack{P \mathrm{\ prime }\\\deg P \le m}} \sum_{n=1}^{\infty} \frac{a_P^n}{n |P|^{ns}} \Bigg| < \varepsilon.
\]
\end{lemma}

\begin{proof}
Let $C = C(q, K, U) > 0$ denote the constant from Lemma~\ref{lem_baby_fundamental}. Let $\ell_0 \ge 1$ and $\varepsilon > 0$ be given. We first control the higher-degree and prime-power tail contributions. Choose an integer $\ell > \ell_0$ sufficiently large so that
\begin{equation}\label{eq_ell_condition}
C q^{-\frac{\ell}{2}} < \frac{\varepsilon}{2} \quad \text{and} \quad \sum_{\substack{P \text{ prime}\\\ell < \deg P}} \; \sum_{n=2}^{\infty} \frac{1}{n |P|^{n \sigma_1}} < \frac{\varepsilon}{2}.
\end{equation}

To verify that such a choice of $\ell$ is possible for the above sum, note that for any $s \in K$ we have $\Re(s) > \sigma_1 > \frac12$. Hence, we obtain
\begin{align*}
\sum_{\substack{P \text{ prime} \\ \ell < \deg P }} \sum_{n=2}^{\infty} \frac{1}{n |P|^{n \sigma_1}}
&\le \frac12 \sum_{d=\ell+1}^{\infty} \; \sum_{\substack{P \text{ prime} \\ \deg P = d}}\; \sum_{n=2}^{\infty} q^{-n d \sigma_1} \\
&\le \frac12 \sum_{d=\ell+1}^{\infty} \pi_q(d) \frac{q^{-2d\sigma_1}}{1 - q^{-d\sigma_1}}.
\end{align*}
Since $\ell > \ell_0 \ge 1$, every integer $d \ge \ell + 1$ satisfies $d \ge 2$ and hence $d \sigma_1 > 1$. Using $\pi_q(d) \le \frac{q^d}{d}$ and the fact that $1 - q^{-d\sigma_1} > 1 - q^{-1} \ge \frac23 > \frac12$ for odd $q \ge 3$, we have
\begin{align}\label{eq_m_tail_bound}
\sum_{\substack{P \text{ prime} \\ \ell < \deg P }} \sum_{n=2}^{\infty} \frac{1}{n |P|^{n \sigma_1}}
&\le \sum_{d=\ell+1}^{\infty} \frac{q^d}{d} q^{-2d\sigma_1} \nonumber \\
&\le \frac{1}{\ell+1} \sum_{d=\ell+1}^{\infty} q^{-d(2\sigma_1 - 1)} \nonumber \\
&\le \frac{q^{-\ell(2\sigma_1 - 1)}}{(2\sigma_1 - 1)\log q},
\end{align}
where the final inequality follows from $e^x - 1 > x$ for $x = (2\sigma_1 - 1)\log q > 0$. Clearly this can be made less than $\frac{\varepsilon}{2}$ by choosing $\ell$ sufficiently large. With $\ell$ now fixed, we introduce the auxiliary function
\[
\widetilde{F}(s) = F(s) - \sum_{\substack{P \text{ prime}\\\deg P \le \ell}} \;\sum_{n=1}^{\infty} \frac{1}{n |P|^{ns}}.
\]
Because the outer sum is finite and each inner series converges uniformly on $U$ (since $|P|^{-n\Re(s)} < q^{-n\sigma_1}$ for all $s \in U$), the function $\widetilde{F}(s)$ is holomorphic on $U$ and real-valued on $U \cap \R$. We now apply Lemma~\ref{lem_baby_fundamental} to $\widetilde{F}(s)$ with this choice of $\ell$. Therefore, there exists a constant $m_0 = m_0(q, U, K, \ell, \widetilde{F}) > \ell > \ell_0$ such that for any integer $m \ge m_0$, there exist numbers $a_P \in \{1,-1\}$ for each prime $P$ with $\ell < \deg P \le m$ satisfying
\begin{equation}\label{eq_cp_approx}
\max_{s \in K} \Bigg| \widetilde{F}(s) - \sum_{\substack{P \text{ prime}\\\ell < \deg P \le m}} \frac{a_P}{|P|^s} \Bigg| \le C q^{-\frac{\ell}{2}} < \frac{\varepsilon}{2}.
\end{equation}
The last inequality follows from~\eqref{eq_ell_condition}. Fixing an integer $m \ge m_0$ and setting $a_P = 1$ for every prime $P$ satisfying $\deg P \le \ell$, we find
\begin{align*}
&\max_{s \in K} \Bigg| F(s) - \sum_{\substack{P \text{ prime}\\\deg P \le m}} \sum_{n=1}^{\infty} \frac{a_P^n}{n |P|^{ns}} \Bigg| \\
&\quad \le \max_{s \in K} \Bigg| \bigg(F(s) - \sum_{\substack{P \text{ prime}\\\deg P \le \ell}} \sum_{n=1}^{\infty} \frac{1}{n |P|^{ns}} \bigg) - \sum_{\substack{P \text{ prime}\\\ell < \deg P \le m}} \frac{a_P}{|P|^s} \Bigg| \\
&\qquad + \max_{s \in K} \Bigg|\sum_{\substack{P \text{ prime}\\\ell < \deg P \le m}} \sum_{n=2}^{\infty} \frac{a_P^n}{n |P|^{ns}}\Bigg| \\
&\quad \le \max_{s \in K} \Bigg| \widetilde{F}(s) - \sum_{\substack{P \text{ prime}\\\ell < \deg P \le m}} \frac{a_P}{|P|^s} \Bigg|
+ \sum_{\substack{P \text{ prime}\\\ell < \deg P}} \;\sum_{n=2}^{\infty} \frac{1}{n |P|^{n \sigma_1}} \\
&\quad < \varepsilon.
\end{align*}
Here, the first inequality holds by the triangle inequality together with the fact that $|a_P| = 1$ and $\Re(s) > \sigma_1$ for all $s \in K$. Combining~\eqref{eq_ell_condition} and~\eqref{eq_cp_approx} establishes the final bound, which completes the proof.
\end{proof}

\section{Results on Quadratic Characters \texorpdfstring{$\chi_D$}{chiD}}\label{Sec quad}

In this section our goal is to prove Lemma~\ref{lem_Agm}, which is essential for our proof of Theorem~\ref{thm 1}. First we record two function field results that we use throughout this section.

\begin{lemma}\label{lem_ensemble_size}
Let $\Hg$ be the set of monic, square-free polynomials $D \in \Fqx$ of degree $2g+1$. Then for any integer $g \ge 1$, the number of elements of the hyperelliptic ensemble is
\[
    \# \Hg = q^{2g}(q - 1).
\]
\end{lemma}
For a proof, see Proposition~2.3 of~\cite{Rosen}.

The next lemma is Lemma~6.2 of~\cite{AK}. We include a proof for the sake of completeness.

\begin{lemma}\label{lem_character_sum}
Let $\chi$ be a non-trivial Dirichlet character modulo $Q$. Then for $n < \deg Q$,
\[
\Bigg| \sum_{\substack{M \ \mathrm{monic}\\\deg M = n}} \chi(M) \Bigg| \le \binom{\deg Q - 1}{n} q^{\frac{n}{2}},
\]
where the sum is over all monic polynomials $M$ of degree $n$, and for $n \ge \deg Q$, the character sum vanishes.
\end{lemma}

\begin{proof}
For $n \ge \deg Q$, write each monic polynomial $M$ of degree $n$ uniquely as $M = A_0 Q + R$, where $A_0$ is monic of degree $n - \deg Q$ and $\deg R < \deg Q$. As $M$ runs through all monic polynomials of degree $n$, the quotient $A_0$ and the remainder $R$ run independently through all $q^{n - \deg Q}$ monic polynomials of degree $n - \deg Q$ and a complete set of residue classes modulo $Q$, respectively. Because $\chi$ is a non-trivial character modulo $Q$, the sum over all residue classes vanishes; that is,
\[
\sum_{\substack{M \text{ monic}\\ \deg M = n}} \chi(M) = q^{n-\deg Q}\sum_{R\pmod Q} \chi(R) = 0.
\]
Hence, the character sum vanishes for all $n \ge \deg Q$.

For $n < \deg Q$, we consider the $u$-form of the general Dirichlet $L$-function in~\eqref{eq_lsdef} with $q^{-s}$ replaced by $u$, namely
\[
\mathcal{L}(u, \chi) = \sum_{M \text{ monic}} \chi(M) u^{\deg M} = \sum_{r=0}^{d} c_r u^r,
\]
where $\displaystyle c_r = \sum_{\substack{M\ \text{monic}\\\deg M = r}} \chi(M)$ and $c_0 = 1$. Because $c_r = 0$ for $r \ge \deg Q$, $\mathcal{L}(u, \chi)$ is a polynomial of degree $d \le \deg Q - 1$. Over $\mathbb{C}$, this polynomial factors as
\[
\mathcal{L}(u, \chi) = \prod_{j=1}^{d} (1 - \alpha_j u).
\]
By the Riemann Hypothesis for curves over finite fields (see Weil's theorem~\cite{Aweil}), the reciprocals of the roots satisfy $|\alpha_j| \le q^{\frac12}$ for all $1 \le j \le d$. If $n > d$, then $c_n = 0$ and the desired bound holds trivially. For $n \le d$, expanding the product allows us to express the coefficient $c_n$ as
\[
c_n = (-1)^n \sum_{1 \le j_1 < \dots < j_n \le d} \alpha_{j_1} \cdots \alpha_{j_n}.
\]
Therefore,
\[
|c_n| \le \binom{d}{n} q^{\frac{n}{2}} \le \binom{\deg Q - 1}{n} q^{\frac{n}{2}}
\]
for every $n < \deg Q$. This completes the proof.
\end{proof}

\begin{lemma}\label{lem_prescribed_characters}
Let $m \ge 1$ be a fixed integer and let a sequence of numbers $a_P \in \{-1, 1\}$ be specified for each prime $P$ with $\deg P \le m$. Define the sub-ensemble
\[
\Hg^m = \left\{ D \in \Hg : \chi_D(P) = a_P \text{ for all primes } P \text{ with } \deg P \le m \right\}.
\]
Then, as $g \to \infty$,
\begin{equation}\label{eq_chidp_ap}
    \frac{\# \Hg^m}{\# \Hg} = C_m + O\left( q^{-\left(g - \frac12\right)} \right),
\end{equation}
where the implied constant in the $O$-term depends at most on $m$ and $q$, and the constant $C_m$ is given by
\begin{equation}\label{eq_Cm}
C_m = \prod_{\substack{P\ \mathrm{prime} \\ \deg P \le m}} \frac12\left(1 + \frac{1}{|P|}\right)^{-1}.
\end{equation}
\end{lemma}

\begin{proof}
For any $D \in \Hg$ and monic prime $P$, recall that the quadratic residue symbol is defined by
\[
\chi_D(P) = \Bigl( \dfrac{D}{P} \Bigr) =
    \begin{cases}
        0,  & \text{if } P \mid D, \\
        1,  & \text{if } P \nmid D \text{ and } D \text{ is a square mod } P, \\
        -1, & \text{if } P \nmid D \text{ and } D \text{ is a non-square mod } P.
    \end{cases}
\]

Let $\mathcal{Q}_P$ and $\mathcal{Q}'_P$ denote the sets of nonzero quadratic residue and non-residue classes modulo $P$, respectively. Because $(\Fqx/P\Fqx)^\times$ is cyclic of even order $|P| - 1$, we have
\[
\# \mathcal{Q}_P = \# \mathcal{Q}'_P = \frac{|P| - 1}{2}.
\]
Consequently, when $a_P = 1$, the condition $\chi_D(P) = a_P$ is satisfied if and only if $D \pmod P \in \mathcal{Q}_P$, and when $a_P = -1$, the condition $\chi_D(P) = a_P$ is satisfied if and only if $D \pmod P \in \mathcal{Q}'_P$.

For our given integer $m \ge 1$, we define the monic square-free modulus
\begin{equation}\label{def Q}
Q = \prod_{\substack{P\ \text{prime} \\ \deg P \le m}} P.
\end{equation}
Suppose now that we have a specific assignment of values $\pm 1$ for $a_P$ for each monic prime $P$ with $\deg P \le m$. For each such $P$, there are $(|P|-1)/2$ residue classes modulo $P$ for which $\chi_D(P) = a_P$, namely the classes in $\mathcal{Q}_P$ if $a_P = 1$ and the classes in $\mathcal{Q}'_P$ if $a_P = -1$. Thus, if $\mathcal{R}$ denotes the set of all residue classes modulo $Q$ for which we simultaneously have $\chi_D(P) = a_P$ for all monic primes $P$ with $\deg P \le m$, then by the Chinese Remainder Theorem we see that
\begin{equation}\label{def_R}
\# \mathcal{R} = \prod_{\substack{P\ \text{prime} \\ \deg P \le m}} \frac{|P| - 1}{2}.
\end{equation}
Since the zero residue class modulo $P$ is excluded from $\mathcal{Q}_P$ and $\mathcal{Q}'_P$ for every $P \mid Q$, every residue class $[C] \in \mathcal{R}$ satisfies $(C, Q) = 1$. Thus, we have
\begin{equation}\label{eq_R_union_split}
\Hg^m = \bigcup_{[C] \in \mathcal{R}} \left\{ D \in \Hg : D \equiv C \pmod Q \right\}
\end{equation}
and hence
\begin{equation}\label{eq_R_sum_split}
\# \Hg^m = \sum_{[C] \in \mathcal{R}} \sum_{\substack{D \in \Hg \\ D \equiv C \pmod Q}} 1.
\end{equation}

Assume that $g$ is sufficiently large so that $2g + 1 \ge \deg Q$. For a fixed class $[C] \in \mathcal{R}$, the orthogonality of Dirichlet characters modulo $Q$ implies that
\[
    \sum_{\substack{D \in \Hg \\ D \equiv C \pmod Q}} 1
    = \sum_{\substack{D\ \text{monic}\\ \deg D=2g+1 \\ D \equiv C \pmod Q}} \mu(D)^2
    = \frac{1}{\Phi(Q)} \sum_{\chi \pmod Q} \overline{\chi}(C) \sum_{\substack{D\ \text{monic}\\ \deg D=2g+1}} \chi(D) \mu(D)^2,
\]
where $\Phi(Q) = |Q| \prod_{P \mid Q} \left(1 - \frac{1}{|P|}\right)$ is the Euler totient function on $\Fqx$. Using the identity $\mu(D)^2 = \sum_{N^2 \mid D} \mu(N)$ and factoring $D = N^2 M$ with $N$ and $M$ monic, $\deg D = 2\deg N + \deg M = 2g+1$, and $(D, Q) = 1$, we find that the inner sum becomes
\[
    \sum_{\substack{D\ \text{monic}\\ \deg D=2g+1}} \chi(D) \mu(D)^2
    = \sum_{\substack{N\ \text{monic}\\\deg N \le g \\ (N, Q) = 1}} \mu(N) \chi(N^2) \sum_{\substack{M\ \text{monic}\\ \deg M=2g+1-2\deg N \\ (M, Q) = 1}} \chi(M).
\]
Therefore,
\begin{equation}\label{eq_orthogonality 1}
     \sum_{\substack{D \in \Hg \\ D \equiv C \pmod Q}} 1
     = \frac{1}{\Phi(Q)} \sum_{\chi \pmod Q} \overline{\chi}(C)\sum_{\substack{N\ \text{monic}\\\deg N \le g \\ (N, Q) = 1}} \mu(N) \chi(N^2) \sum_{\substack{M\ \text{monic}\\ \deg M=2g+1-2\deg N \\ (M, Q) = 1}} \chi(M).
\end{equation}
We denote the contribution from the principal character as $S_1$, and the contribution from the other characters as $S_2$, so that
\begin{equation}\label{eq_orthogonality 2}
 \sum_{\substack{D \in \Hg \\ D \equiv C \pmod Q}} 1 = S_1 + S_2.
\end{equation}

From~\eqref{eq_orthogonality 1}, the contribution from the principal character is therefore
\[
    S_1 = \frac{1}{\Phi(Q)} \sum_{\substack{N\ \text{monic}\\\deg N \le g \\ (N, Q) = 1}} \mu(N) \sum_{\substack{M\ \text{monic}\\ \deg M=2g+1-2\deg N \\ (M, Q) = 1}} 1.
\]
To evaluate the inner sum over $M$, set $d = 2g + 1 - 2\deg N$ and apply the Möbius inversion identity $\sum_{E \mid (M, Q)} \mu(E) = 1$ if $(M, Q) = 1$ and $0$ otherwise (where $E$ runs over monic polynomials). Interchanging the order of summation gives
\begin{equation}\label{eq_S1_mobius}
    \sum_{\substack{M\ \text{monic}\\\deg M = d \\ (M, Q) = 1}} 1
    = \sum_{\substack{M\ \text{monic}\\\deg M = d}} \quad \sum_{\substack{E\ \text{monic} \\ E \mid (M,Q)}} \mu(E)
    = \sum_{\substack{E\ \text{monic} \\ E \mid Q}} \mu(E) \sum_{\substack{M\ \text{monic} \\ \deg M = d \\ E \mid M}} 1.
\end{equation}
A monic polynomial $M$ of degree $d$ is divisible by $E$ if and only if $M = E B$ for a unique monic polynomial $B$ of degree $d - \deg E$. Because no nonzero polynomial has negative degree, the innermost sum vanishes whenever $\deg E > d$, so that
\[
    \sum_{\substack{M\ \text{monic} \\ \deg M = d \\ E \mid M}} 1 =
    \begin{cases}
        q^{d - \deg E}, & \text{if } \deg E \le d, \\[4pt]
        0, & \text{if } \deg E > d.
    \end{cases}
\]
By~\eqref{def Q}, the modulus $Q$ is square-free, so its largest monic divisor is $E = Q$ itself. Consequently, the condition $\deg E \le d$ holds for every monic divisor $E \mid Q$ if and only if $d \ge \deg Q$, which is equivalent to
\[
    \deg N \le g + \frac12 - \frac12 \deg Q.
\]
When $\deg N$ exceeds this threshold, divisors $E \mid Q$ with $\deg E > d$ contribute $0$ rather than $q^{d - \deg E}$ in~\eqref{eq_S1_mobius}, preventing the sum from factoring cleanly into $q^d \Phi(Q)/|Q|$.

We therefore split the sum over $N$ into two parts, writing $S_1 = S_{1,1} + S_{1,2}$, where $S_{1,1}$ contains the terms with $\deg N \le g + \frac12 - \frac12 \deg Q$ and $S_{1,2}$ contains the terms with $g + \frac12 - \frac12 \deg Q < \deg N \le g$. For every $N$ in $S_{1,1}$, we have $\deg E \le \deg Q \le d$ for all monic divisors $E \mid Q$, so~\eqref{eq_S1_mobius} evaluates exactly as
\begin{align*}
    \sum_{\substack{M\ \text{monic}\\\deg M=2g+1-2\deg N \\ (M, Q) = 1}} 1
    &= \sum_{\substack{E\ \text{monic} \\ E \mid Q}} \mu(E)\; q^{2g+1 - 2\deg N-\deg E} \\
    &= q^{2g+1 - 2\deg N} \sum_{\substack{E\ \text{monic} \\ E \mid Q}} \frac{\mu(E)}{|E|} \\
    &= q^{2g+1 - 2\deg N}\; \frac{\Phi(Q)}{|Q|}.
\end{align*}
Substituting this identity into $S_{1,1}$ and canceling the factor $\Phi(Q)$, we obtain
\begin{align*}
    S_{1,1}
    &= \frac{q^{2g+1}}{|Q|} \sum_{\substack{N\ \text{monic} \\ \deg N \le g + \frac12 - \frac12 \deg Q \\ (N, Q) = 1}} \frac{\mu(N)}{|N|^2} \\
    &= \frac{q^{2g+1}}{|Q|} \sum_{\substack{N\ \text{monic} \\ (N, Q) = 1}} \frac{\mu(N)}{|N|^2} + O\Bigg( \frac{q^{2g+1}}{|Q|} \sum_{\substack{N\ \text{monic} \\ \deg N > g + \frac12 - \frac12 \deg Q}} \frac{1}{|N|^2} \Bigg).
\end{align*}
Writing the first series as an Euler product over the primes $P$, we find that
\[
    \sum_{\substack{N\ \text{monic} \\ (N, Q) = 1}} \frac{\mu(N)}{|N|^2}
    = \zeta_{\Fqx}(2)^{-1} \prod_{\substack{P\ \text{prime} \\ P \mid Q}} \left(1 - \frac{1}{|P|^2}\right)^{-1}.
\]
By~\eqref{eq_zeta_closed}, we have $\zeta_{\Fqx}(2) = \frac{q}{q-1}$. On the other hand, since $g + \frac12 - \frac12 \deg Q \ge 0$ and $q \ge 3$, the tail series in the $O$-term is bounded by
\[
    \frac{q^{2g+1}}{|Q|} \sum_{k > g + \frac12 - \frac12 \deg Q} \frac{1}{q^{k}}
    \le \frac{q^{2g+1}}{|Q|} \frac{q}{q-1} q^{-\left(g + \frac12 - \frac12 \deg Q\right)}
    \le \frac{3}{2} \frac{q^{g + \frac12}}{|Q|^{\frac12}}
    \le \frac{3}{2} |Q| q^{g + \frac12}.
\]
For $S_{1,2}$, dropping the coprimality condition $(M, Q) = 1$ and bounding the inner sum over $M$ trivially by $q^{2g+1-2\deg N}$, we obtain
\begin{align*}
    |S_{1,2}|
    &\le \frac{1}{\Phi(Q)} \sum_{\substack{N\ \text{monic} \\ g + \frac12 - \frac12 \deg Q < \deg N \le g}} q^{2g+1-2\deg N} \\
    &= \frac{q^{2g+1}}{\Phi(Q)} \sum_{g + \frac12 - \frac12 \deg Q < j \le g} q^{-j}
    \le \frac{3}{2} \frac{|Q|^{\frac12}}{\Phi(Q)} q^{g + \frac12}
    \le \frac{3}{2} |Q| q^{g + \frac12}.
\end{align*}
Combining $S_{1,1}$ and $S_{1,2}$, we conclude that
\[
    S_1 = \frac{q^{2g}(q-1)}{|Q|}\prod_{\substack{P\ \text{prime} \\ P \mid Q}} \left(1 - \frac{1}{|P|^2}\right)^{-1} + O\left( |Q| q^{g + \frac12} \right).
\]

To treat $S_2$, note that by Lemma~\ref{lem_character_sum}, the character sum over monic polynomials $M$ vanishes whenever the degree of $M$ is at least $\deg Q$. Thus
\[
    \sum_{\substack{M\ \text{monic} \\ \deg M=2g+1-2\deg N}} \chi(M) = 0
\]
when $2g + 1 - 2\deg N \ge \deg Q$. Since by~\eqref{eq_orthogonality 1} we also have $\deg N \le g$, nonzero contributions occur only when
\[
    g + \tfrac{1}{2} - \tfrac{1}{2} \deg Q < \deg N \le g.
\]
If $\deg N$ is in this range, Lemma~\ref{lem_character_sum} gives
\[
    \Bigg| \sum_{\substack{M\ \text{monic} \\ \deg M=2g+1-2\deg N}} \chi(M) \Bigg| \le \binom{\deg Q - 1}{2g+1-2\deg N} q^{g - \deg N + \frac12}.
\]
Summing over the $\Phi(Q) - 1$ non-principal characters and grouping the monic polynomials $N$ by their degree, we see from~\eqref{eq_orthogonality 1} that
\begin{align*}
    |S_2|
    &\le \frac{\Phi(Q) - 1}{\Phi(Q)} \sum_{\substack{N\ \text{monic} \\ g + \frac12 - \frac{\deg Q}{2} < \deg N \le g}} \binom{\deg Q - 1}{2g+1-2\deg N} q^{g + \frac12 - \deg N} \\
    &\le q^{g + \frac12} \sum_{g + \frac12 - \frac{\deg Q}{2} < j \le g} \binom{\deg Q - 1}{2g+1-2j} \\
    &\le 2^{\deg Q - 1} q^{g + \frac12} \le |Q| q^{g + \frac12},
\end{align*}
since $2 < q$.

Combining $S_1$ and $S_2$ via~\eqref{eq_orthogonality 2}, we find that for every $[C] \in \mathcal{R}$ and all $g$ with $2g + 1 \ge \deg Q$,
\begin{equation}\label{eq_num_of_coprime_in_hg}
    \sum_{\substack{D \in \Hg \\ D \equiv C \pmod Q}} 1 = \frac{q^{2g}(q-1)}{|Q|} \prod_{\substack{P\ \text{prime} \\ P \mid Q}} \left(1 - \frac{1}{|P|^2}\right)^{-1} + O\left(|Q| q^{g + \frac12} \right),
\end{equation}
where the implied constant is absolute. Summing over all $\#\mathcal{R}$ classes $[C] \in \mathcal{R}$ and using~\eqref{def Q}, \eqref{def_R}, and~\eqref{eq_R_sum_split}, we find that
\begin{align*}
    \# \Hg^m
    &= \# \mathcal{R} \Bigg( \frac{q^{2g}(q-1)}{|Q|} \prod_{\substack{P\ \text{prime} \\ P \mid Q}} \left(1 - \frac{1}{|P|^2}\right)^{-1} + O\Big(|Q| q^{g + \frac12} \Big) \Bigg) \\
    &= q^{2g}(q-1) \prod_{\substack{P\ \text{prime} \\ \deg P \le m}} \Bigg( \frac{|P| - 1}{2|P|} \left(1 - \frac{1}{|P|^2}\right)^{-1} \Bigg) + O\Big( \#\mathcal{R} \cdot |Q| q^{g + \frac12} \Big) \\
    &= q^{2g}(q-1) C_m + O\Big( \#\mathcal{R} \cdot |Q| q^{g + \frac12} \Big).
\end{align*}
Dividing by $\# \Hg = q^{2g}(q-1)$ from Lemma~\ref{lem_ensemble_size} gives
\[
    \frac{\# \Hg^m}{\# \Hg} = C_m + O\left( \frac{\#\mathcal{R} \cdot |Q|}{q - 1} q^{-\left(g - \frac12\right)} \right).
\]
Because both $\#\mathcal{R}$ and $|Q|$ depend only on $m$ and $q$, this establishes~\eqref{eq_chidp_ap} with an implied constant depending at most on $m$ and $q$. This completes the proof.
\end{proof}

\begin{lemma}\label{lem_mid_bound}
Let $m \ge 1$ be a fixed integer and let a sequence of numbers $a_P \in \{-1, 1\}$ be specified for each prime $P$ with $\deg P \le m$. Let $\Hg^m$ and $C_m$ be defined as in Lemma~\ref{lem_prescribed_characters}. Let $h_g \ge 1$ be an integer satisfying $h_g \le \alpha g$ with $0 < \alpha \le \frac{1}{2(1 + 2\log_q 2)}$, and fix $\sigma_1 > \frac12$. Then there exists an integer $g_0 = g_0(m, q, \sigma_1) \ge 1$ such that for any $g \ge g_0$ and uniformly for $s \in \mathbb{C}$ with $\Re(s) \ge \sigma_1$,
\[
\sum_{D \in \Hg^m} \Bigg| \sum_{\substack{P\ \mathrm{prime} \\ m < \deg P \le h_g}} \frac{\chi_D(P)}{|P|^s} \Bigg|^2 \ll \frac{(q^m)^{1 - 2\sigma_1}}{(2\sigma_1 - 1)\log(q^m)} \, C_m \, \#\Hg,
\]
where the implied constant is absolute.
\end{lemma}

\begin{proof}
For a given $m \ge 1$, let $Q$ and $\mathcal{R}$ be as in~\eqref{def Q} and~\eqref{def_R}, and assume $g$ is large enough that $2g + 1 \ge \deg Q$. Then by~\eqref{eq_R_union_split} we have
\begin{equation}\label{eq_mid_sum}
\sum_{D \in \Hg^m} \Bigg| \sum_{\substack{P\ \text{prime}\\m < \deg P \le h_g}} \frac{\chi_D(P)}{|P|^s} \Bigg|^2
= \sum_{[C] \in \mathcal{R}} \; \sum_{\substack{D \in \Hg \\ D \equiv C \pmod{Q}}} \Bigg| \sum_{\substack{P\ \text{prime}\\m < \deg P \le h_g}} \frac{\chi_D(P)}{|P|^s} \Bigg|^2.
\end{equation}

Setting $\sigma = \Re(s)$ and expanding the squared modulus inside the inner sum, we find that
\begin{align}\label{eq_S1C+S2C}
\sum_{\substack{D \in \Hg \\ D \equiv C \pmod{Q}}} \Bigg| \sum_{\substack{P\ \text{prime}\\m < \deg P \le h_g}} \frac{\chi_D(P)}{|P|^s} \Bigg|^2
&= \sum_{\substack{D \in \Hg \\ D \equiv C \pmod{Q}}} \;\; \sum_{\substack{P_1, P_2\ \text{prime}\\m < \deg P_1,\, \deg P_2 \le h_g}} \frac{\chi_D(P_1)\overline{\chi_D(P_2)}}{|P_1|^s |P_2|^{\bar{s}}} \nonumber \\
&= \sum_{\substack{P\ \text{prime}\\m < \deg P \le h_g}} \frac{1}{|P|^{2\sigma}} \bigg(\sum_{\substack{D \in \Hg \\ D \equiv C \pmod{Q}}} |\chi_D(P)|^2\bigg) \nonumber \\
&\quad + \sum_{\substack{P_1, P_2\ \text{prime}\\P_1 \neq P_2 \\ m < \deg P_1,\, \deg P_2 \le h_g}} \frac{1}{|P_1|^s |P_2|^{\bar{s}}} \bigg( \sum_{\substack{D \in \Hg \\ D \equiv C \pmod{Q}}} \chi_D(P_1)\overline{\chi_D(P_2)}\bigg) \nonumber \\
&= S_1(C) + S_2(C),
\end{align}

where $S_1(C)$ denotes the diagonal contribution and $S_2(C)$ denotes the off-diagonal contribution.

For the diagonal sum $S_1(C)$, we see that
\[
|S_1(C)| \le \bigg(\sum_{\substack{P\ \text{prime}\\m < \deg P \le h_g}} \frac{1}{|P|^{2\sigma}}\bigg) \bigg( \sum_{\substack{D \in \Hg \\ D \equiv C \pmod{Q}}} 1 \bigg).
\]
For $\sigma \ge \sigma_1 > \frac12$, using the exact upper bound $\pi_q(j) \le \frac{q^j}{j}$, we obtain
\begin{align*}
\sum_{\substack{P\ \text{prime}\\m < \deg P \le h_g}} \frac{1}{|P|^{2\sigma}}
&\le \sum_{m < j \le h_g} \frac{q^j}{j \cdot q^{2\sigma_1 j}} \\
&\le \frac{1}{m} \int_{m}^{\infty} q^{-u(2\sigma_1 - 1)} \, du \\
&= \frac{(q^m)^{1 - 2\sigma_1}}{(2\sigma_1 - 1)\log(q^m)}.
\end{align*}
From this and~\eqref{eq_num_of_coprime_in_hg}, we obtain
\begin{equation}\label{eq_S1C}
S_1(C) \ll \frac{(q^m)^{1 - 2\sigma_1}}{(2\sigma_1 - 1)\log(q^m)} \Bigg( \frac{q^{2g}(q-1)}{|Q|} \prod_{\substack{P\ \text{prime} \\ P \mid Q}} \left(1 - \frac{1}{|P|^2}\right)^{-1} + |Q| q^{g + \frac12} \Bigg),
\end{equation}
with an absolute implied constant.

We now bound the sum $S_2(C)$. For fixed distinct primes $P_1 \neq P_2$ with $m < \deg P_1, \deg P_2 \le h_g$, the orthogonality of Dirichlet characters modulo $Q$ gives
\begin{equation}\label{eq_sum_chi_chi}
\sum_{\substack{D \in \Hg \\ D \equiv C \pmod{Q}}} \chi_D(P_1)\overline{\chi_D(P_2)}
= \frac{1}{\Phi(Q)} \sum_{\chi \pmod{Q}} \overline{\chi}(C)\sum_{\substack{D\ \text{monic}\\ \deg D=2g+1}} \mu(D)^2 \legendre{D}{P_1 P_2} \chi(D).
\end{equation}
Using the identity $\displaystyle \mu(D)^2 = \sum_{\substack{N\ \text{monic} \\ N^2 \mid D}} \mu(N)$ and factoring $D = N^2 M$ with $N$ and $M$ monic, we decompose the innermost sum as
\begin{align*}
\sum_{\substack{D\ \text{monic}\\ \deg D=2g+1}} \mu(D)^2 &\left(\frac{D}{P_1 P_2}\right) \chi(D) \\
&= \sum_{\substack{N\ \text{monic}\\\deg N \le g }} \mu(N) \left(\frac{N^2}{P_1 P_2}\right) \chi(N^2)
\Bigg( \sum_{\substack{M\ \text{monic}\\ \deg M = 2g+1 - 2\deg N}} \left(\frac{M}{P_1 P_2}\right) \chi(M)\Bigg).
\end{align*}
Because $\deg P_1, \deg P_2 > m$, both $P_1$ and $P_2$ are coprime to $Q$, and $P_1 \neq P_2$. Consequently, the character $\psi(M) = \legendre{M}{P_1 P_2}\chi(M)$ is a non-trivial Dirichlet character modulo $Q P_1 P_2$ for every character $\chi \pmod Q$, including the principal character. By Lemma~\ref{lem_character_sum}, the sum over monic polynomials $M$ of degree $n = 2g+1 - 2\deg N$ vanishes identically when $n \ge \deg(Q P_1 P_2)$. Thus, nonzero contributions occur exclusively for monic polynomials $N$ in the degree range
\[
g + \frac12 - \frac{\deg(Q P_1 P_2)}{2} < \deg N \le g.
\]
Therefore, applying Lemma~\ref{lem_character_sum} again and summing over monic polynomials $N$ grouped by degree, we find that
\begin{align*}
\Bigg| \sum_{\substack{D\ \text{monic}\\ \deg D=2g+1}} \mu(D)^2 \left(\frac{D}{P_1 P_2}\right) \chi(D)\Bigg|
&\le \sum_{\substack{N\ \text{monic}\\ g+\frac12 - \frac12\deg(QP_1P_2) < \deg N \le g }}
\binom{\deg(QP_1P_2) - 1}{2g+1-2\deg N} q^{g - \deg N + \frac12} \\
&\le q^{g+\frac12} \sum_{g+\frac12 - \frac12\deg(QP_1P_2) < j \le g}
\binom{\deg(QP_1P_2) - 1}{2g+1-2j} \\
&\le q^{g + \frac12} 2^{\deg(Q P_1 P_2) - 1} \ll q^{g + \frac12} 2^{2h_g} |Q|^{\log_q 2}.
\end{align*}
Using this bound and~\eqref{eq_sum_chi_chi}, we find that
\[
\sum_{\substack{D \in \Hg \\ D \equiv C \pmod{Q}}} \chi_D(P_1)\overline{\chi_D(P_2)}
\ll q^{g + \frac12} 2^{2h_g} |Q|^{\log_q 2}.
\]
Inserting this estimate into~\eqref{eq_S1C+S2C} and recalling that $\sigma \ge \sigma_1 > \frac12$, we obtain
\begin{align*}
|S_2(C)| &\ll q^{g + \frac12} 2^{2h_g} |Q|^{\log_q 2} \Bigg( \sum_{\substack{P\ \text{prime}\\m < \deg P \le h_g}} \frac{1}{|P|^{\frac12}} \Bigg)^2 \\
&\ll q^{g + \frac12} 2^{2h_g} |Q|^{\log_q 2} \Bigg( \sum_{m < j \le h_g} \frac{q^{\frac{j}{2}}}{j} \Bigg)^2 \\
&\ll q^{g + h_g(1 + 2\log_q 2) + \frac12} |Q|^{\log_q 2}.
\end{align*}
In the second step we have used $\pi_q(j) \le \frac{q^j}{j}$. By our hypothesis, $h_g \le \alpha g$ with $0 < \alpha \le \frac{1}{2(1 + 2\log_q 2)}$, so the exponent of $q$ above satisfies
\[
g + h_g(1 + 2\log_q 2) + \frac12 \le g\left(1 + \alpha(1 + 2\log_q 2)\right) + \frac12 \le \frac{3}{2}g + \frac12.
\]
Thus,
\[
|S_2(C)| \ll q^{\frac{3}{2}g + \frac12} |Q|^{\log_q 2}.
\]
Combining this bound and~\eqref{eq_S1C} in~\eqref{eq_S1C+S2C}, we see that for any fixed $m$, $q$, and $\sigma_1$, there exists a sufficiently large integer $g_0 = g_0(m, q, \sigma_1)$ such that whenever $g \ge g_0$, we have
\begin{align*}
\sum_{\substack{D \in \Hg \\ D \equiv C \pmod{Q}}} \Bigg| \sum_{\substack{P\ \text{prime}\\m < \deg P \le h_g}} \frac{\chi_D(P)}{|P|^s} \Bigg|^2
&\ll \frac{(q^m)^{1 - 2\sigma_1}}{(2\sigma_1 - 1)\log(q^m)} \Bigg( \frac{q^{2g}(q-1)}{|Q|} \prod_{\substack{P\ \text{prime} \\ P \mid Q}} \left(1 - \frac{1}{|P|^2}\right)^{-1} + |Q| q^{g + \frac12}\Bigg) \\
&\qquad + q^{\frac{3}{2}g + \frac12} |Q|^{\log_q 2} \\
&\ll \frac{(q^m)^{1 - 2\sigma_1}}{(2\sigma_1 - 1)\log(q^m)} \Bigg( \frac{q^{2g}(q-1)}{|Q|} \prod_{\substack{P\ \text{prime} \\ P \mid Q}} \left(1 - \frac{1}{|P|^2}\right)^{-1} \Bigg),
\end{align*}
with an absolute implied constant.

Now, summing over all $[C] \in \mathcal{R}$ on the right-hand side of~\eqref{eq_mid_sum} and using the formulas (see~\eqref{def_R} and~\eqref{eq_Cm})
\[
\#\mathcal{R} = \prod_{\substack{P\ \text{prime} \\ \deg P \le m}} \frac{|P| - 1}{2}, \qquad C_m = \prod_{\substack{P\ \text{prime}\\\deg P \le m}} \frac12\left(1 + \frac{1}{|P|}\right)^{-1}
\]
along with the relation $\#\Hg = q^{2g}(q-1)$ from Lemma~\ref{lem_ensemble_size}, we obtain for all $g \ge g_0$ that
\[
\sum_{D \in \Hg^m} \Bigg| \sum_{\substack{P\ \text{prime}\\m < \deg P \le h_g}} \frac{\chi_D(P)}{|P|^s} \Bigg|^2
\ll \frac{(q^m)^{1 - 2\sigma_1}}{(2\sigma_1 - 1)\log(q^m)} \, C_m \, \#\Hg.
\]
This concludes the proof.
\end{proof}

\begin{lemma}\label{lem_Agm}
Let $\sigma_1 > \frac12$ and $\varepsilon > 0$, and let $K \subset \mathbb{C}$ be a non-empty compact set such that $K \subset \{s \in \mathbb{C} : \Re(s) > \sigma_1\}$. Let $h_g \ge 1$ be an integer satisfying $h_g \le \alpha g$ with $0 < \alpha \le \frac{1}{2(1 + 2\log_q 2)}$. Then there exists an integer $m_0 = m_0(q, \sigma_1, K, \varepsilon) \ge 1$ such that for any fixed $m \ge m_0$, the following is true. For any choice of numbers $a_P \in \{1, -1\}$ with $\deg P \le m$, the set
\[
\mathcal{A}_{2g+1}^m = \mathcal{A}_{2g+1}^m(K, \varepsilon, \{a_P\}) = \Bigg\{ D \in \Hg^m : \max_{s \in K} \Bigg| \sum_{\substack{P\ \mathrm{prime} \\ m < \deg P \le h_g}} \frac{\chi_D(P)}{|P|^s} \Bigg| < \varepsilon \Bigg\}
\]
satisfies, for all sufficiently large $g$,
\[
\frac{\#\mathcal{A}_{2g+1}^m}{\#\Hg} > \frac12 C_m.
\]
\end{lemma}

\begin{proof}
Let $\sigma_3 = 1 + \sup \{\Re(s) : s \in K\}$ and $B = 1 + \sup \{ |\Im(s)| : s \in K\}$. Let $W \subset \mathbb{C}$ be the open rectangular region with vertices $\sigma_1 \pm iB$ and $\sigma_3 \pm iB$ so that $K \subset W$. Let $a(W, K) > 0$ be the constant from Lemma~\ref{lem_Titchmarsh} associated with $W$ and $K$, and let $b > 0$ be the absolute implied constant from Lemma~\ref{lem_mid_bound}. We choose $m_0 = m_0(q, \sigma_1, K, \varepsilon) \ge 1$ sufficiently large so that for all $m \ge m_0$ we have
\begin{equation}\label{eq_m0_condition}
\operatorname{Area}(W) \cdot b \, \frac{(q^{m})^{1 - 2\sigma_1}}{(2\sigma_1 - 1)\log(q^{m})} < \frac{\varepsilon^2}{4 a(W, K)^2}.
\end{equation}

Now fix any $m \ge m_0$ and define the set
\[
\widetilde{\mathcal{A}}_{2g+1}^m = \Bigg\{ D \in \Hg^m : \iint_W \Bigg| \sum_{\substack{P\ \text{prime} \\ m < \deg P \le h_g}} \frac{\chi_D(P)}{|P|^{\sigma + it}} \Bigg|^2 d\sigma \, dt < \frac{\varepsilon^2}{a(W, K)^2} \Bigg\}.
\]
By Lemma~\ref{lem_Titchmarsh}, any polynomial $D \in \widetilde{\mathcal{A}}_{2g+1}^m$ satisfies $\displaystyle \max_{s \in K} \Bigg| \sum_{\substack{P\ \text{prime} \\ m < \deg P \le h_g}} \frac{\chi_D(P)}{|P|^s} \Bigg| < \varepsilon$, and therefore
\[
\widetilde{\mathcal{A}}_{2g+1}^m \subseteq \mathcal{A}_{2g+1}^m.
\]
By Lemma~\ref{lem_mid_bound}, for fixed $m$, $q$, and $\sigma_1$, there exists $g_0 = g_0(m, q, \sigma_1) \ge 1$ such that for all $g \ge g_0$,
\begin{align}\label{eq_sum_integral_bound}
\sum_{D \in \Hg^m} \iint_W \Bigg| \sum_{\substack{P\ \text{prime} \\ m < \deg P \le h_g}} \frac{\chi_D(P)}{|P|^{\sigma + it}} \Bigg|^2 d\sigma \, dt
&\le \operatorname{Area}(W) \cdot b \, \frac{(q^m)^{1 - 2\sigma_1}}{(2\sigma_1 - 1)\log(q^m)} \, C_m \, \#\Hg \nonumber \\
&< \frac{\varepsilon^2}{4 a(W, K)^2} \, C_m \, \#\Hg.
\end{align}
The last step follows from~\eqref{eq_m0_condition}.

We now wish to bound the cardinality of the difference set $\Hg^m \setminus \widetilde{\mathcal{A}}_{2g+1}^m$. From the definition of $\widetilde{\mathcal{A}}_{2g+1}^m$ we see that
\begin{align*}
\sum_{D \in \Hg^m} \iint_W \Bigg| \sum_{\substack{P\ \text{prime} \\ m < \deg P \le h_g}} \frac{\chi_D(P)}{|P|^{\sigma + it}} \Bigg|^2 d\sigma \, dt
&\ge \sum_{D \in \Hg^m \setminus \widetilde{\mathcal{A}}_{2g+1}^m} \iint_W \Bigg| \sum_{\substack{P\ \text{prime} \\ m < \deg P \le h_g}} \frac{\chi_D(P)}{|P|^{\sigma + it}} \Bigg|^2 d\sigma \, dt \\
&\ge \left( \# \Hg^m - \#\widetilde{\mathcal{A}}_{2g+1}^m \right) \cdot \frac{\varepsilon^2}{a(W, K)^2}.
\end{align*}
Hence, combining this lower bound with inequality~\eqref{eq_sum_integral_bound}, we obtain
\[
\# \Hg^m - \#\widetilde{\mathcal{A}}_{2g+1}^m < \frac{1}{4} C_m \#\Hg.
\]
From this it follows that
\[
\frac{\#\widetilde{\mathcal{A}}_{2g+1}^m}{\#\Hg} > \frac{\#\Hg^m}{\#\Hg} - \frac{1}{4} C_m.
\]
Now, using the inclusion $\widetilde{\mathcal{A}}_{2g+1}^m \subseteq \mathcal{A}_{2g+1}^m$ together with Lemma~\ref{lem_prescribed_characters}, we find that for all sufficiently large $g$,
\[
\frac{\#\mathcal{A}_{2g+1}^m}{\#\Hg} > \frac{3}{4} C_m - \frac{1}{4} C_m = \frac12 C_m.
\]
This completes the proof.
\end{proof}

\section{Proof of Theorem~\ref{thm 1}}

Our goal is to show that for a positive proportion of $D \in \Hg$, the $L$-function $L(s, \chi_D)$ uniformly approximates $F(s)$ on $K$, where $U$, $K$, and $F(s)$ satisfy all the hypotheses of Theorem~\ref{thm 1}. To apply the Fundamental Lemma at the logarithmic level, we first construct a holomorphic logarithm $f(s) = \log F(s)$ on $U$ that is real-valued on $U \cap \R$.

Because $U$ is an open rectangular region symmetric about the real axis, it is connected, simply connected, with $U \cap \R = (\sigma_1, \sigma_2) \neq \emptyset$. Since $F$ is holomorphic and non-vanishing on $U$, Rudin~\cite[Theorem~13.11(h)]{WRudin} provides a function $f_0$ that is holomorphic on $U$ such that $\exp(f_0(s)) = F(s)$ for all $s \in U$. Furthermore, because $F$ is holomorphic on $U$ and real-valued on $U \cap \R$, the Schwarz reflection principle gives $F(\overline{s}) = \overline{F(s)}$ for all $s \in U$. Consequently,
\[
\exp\bigl(\overline{f_0(\overline{s})}\bigr) = \overline{F(\overline{s})} = F(s) = \exp\bigl(f_0(s)\bigr) \qquad (s \in U),
\]
so $f_0(s) - \overline{f_0(\overline{s})} = 2\pi i k$ on the connected domain $U$ for some constant integer $k$. Define
\[
f(s) = \frac{f_0(s) + \overline{f_0(\overline{s})}}{2} = f_0(s) - \pi i k \qquad (s \in U).
\]
Then $f$ is holomorphic on $U$, real-valued on $U \cap \R$, and satisfies $\exp(f(s)) = (-1)^k F(s)$ throughout $U$. Choosing any point $x_0 \in U \cap \R$, we have $f(x_0) \in \R$, so $\exp(f(x_0)) > 0$; since $F(x_0) > 0$, it follows that $(-1)^k = 1$. Therefore, $f$ is holomorphic on $U$, real-valued on $U \cap \R$, and satisfies $\exp(f(s)) = F(s)$ for all $s \in U$.

Next, let $\displaystyle M = \max_{s \in K} |F(s)| > 0$ and choose $\varepsilon_1 > 0$ so small that
\begin{equation}\label{eq_choice_of_epsilon1}
\varepsilon_1 < \min\Big\{ \frac{1}{8}, \, \frac{\varepsilon}{8M + 1} \Big\}.
\end{equation}
For every integer $g \ge 2$, set
\[
h_g = \left\lfloor \frac{2g}{(2\sigma_1 - 1)\log_q g} \right\rfloor.
\]
Since $h_g \to \infty$ and $h_g = o(g)$ as $g \to \infty$, we have $1 \le h_g \le \alpha g$ with $\alpha = \frac{1}{2(1 + 2\log_q 2)}$ for all sufficiently large $g$. Let $m_0 = m_0(q, \sigma_1, K, \varepsilon_1) \ge 1$ be the constant from Lemma~\ref{lem_Agm}, and fix an integer $\ell_0 \ge m_0$ sufficiently large so that for all $m > \ell_0$, we have
\[
\frac{q^{-m(2\sigma_1 - 1)}}{(2\sigma_1 - 1) \log q} < \varepsilon_1.
\]

Now we apply Lemma~\ref{lem_fundamental}, the Fundamental Lemma, to $f(s)$ with $\varepsilon = \varepsilon_1$ and $\ell_0$ as above. We obtain a fixed integer $m = m(q, U, K, \varepsilon_1, f, \ell_0) > \ell_0 \ge m_0$ and a choice of numbers $a_P \in \{1, -1\}$ for each prime $P$ with $\deg P \le m$, all independent of $g$, such that
\begin{equation}\label{eq_fund_lemma_approx}
\max_{s \in K} \Bigg| f(s) - \sum_{\substack{P\ \text{prime} \\ \deg P \le m}} \sum_{n=1}^{\infty} \frac{a_P^n}{n |P|^{ns}} \Bigg| < \varepsilon_1.
\end{equation}

We now restrict our attention to the sub-ensemble of $\Hg$ where $\chi_D(P)$ matches $a_P$ for the primes $P$ up to degree $m$. Recall that
\[
\Hg^m = \left\{ D \in \Hg : \chi_D(P) = a_P \text{ for all primes } P \text{ with } \deg P \le m \right\}.
\]
By Lemma~\ref{lem_Agm}, there exists a strictly positive constant $C_m > 0$ (independent of $g$) such that the subset $\mathcal{A}_{2g+1}^m = \mathcal{A}_{2g+1}^m(K, \varepsilon_1, \{a_P\}) \subseteq \Hg^m$ satisfies, for all sufficiently large $g$,
\begin{equation}\label{eq_positive_density}
\frac{\# \mathcal{A}_{2g+1}^m}{\#\Hg} > \frac12 C_m > 0.
\end{equation}

We next show that for all sufficiently large $g$ (in particular, with $h_g \ge m$) and every $D \in \mathcal{A}_{2g+1}^m$, replacing the series in~\eqref{eq_fund_lemma_approx} by the truncated prime-power sum up to degree $h_g$ introduces an error of at most $3\varepsilon_1$. Since $\chi_D(P) = a_P$ for every prime $P$ with $\deg P \le m$, we have $\chi_D(P^n) = a_P^n$ for all $n \ge 1$ and $\deg P \le m$. Therefore,
\begin{align}\label{eq_m_to_hg}
\sum_{\substack{P\ \text{prime} \\ \deg P \le h_g}} \sum_{n=1}^{\lfloor h_g / \deg P \rfloor} \frac{\chi_D(P^n)}{n |P|^{ns}}
&- \sum_{\substack{P\ \text{prime} \\ \deg P \le m}} \sum_{n=1}^{\infty} \frac{a_P^n}{n |P|^{ns}} \nonumber \\
&= \sum_{\substack{P\ \text{prime} \\ m < \deg P \le h_g}} \frac{\chi_D(P)}{|P|^s}
+ \sum_{\substack{P\ \text{prime} \\ m < \deg P \le h_g}} \sum_{n=2}^{\lfloor h_g / \deg P \rfloor} \frac{\chi_D(P^n)}{n |P|^{ns}}\\
&\qquad - \sum_{\substack{P\ \text{prime} \\ \deg P \le m}} \sum_{n > \lfloor h_g / \deg P \rfloor} \frac{a_P^n}{n |P|^{ns}}.
\end{align}
By the definition of $\mathcal{A}_{2g+1}^m$ in Lemma~\ref{lem_Agm}, every $D \in \mathcal{A}_{2g+1}^m$ satisfies
\[
\max_{s \in K} \Bigg| \sum_{\substack{P\ \text{prime} \\ m < \deg P \le h_g}} \frac{\chi_D(P)}{|P|^s} \Bigg| < \varepsilon_1.
\]
For the second sum on the right-hand side of~\eqref{eq_m_to_hg}, the same argument used in~\eqref{eq_m_tail_bound} (with $\ell$ replaced by $m > \ell_0 \ge 1$) shows that for all $s \in K$,
\[
\sum_{\substack{P\ \text{prime} \\ m < \deg P \le h_g}} \sum_{n=2}^{\lfloor h_g / \deg P \rfloor} \frac{1}{n |P|^{n \sigma_1}} \le \frac{q^{-m(2\sigma_1 - 1)}}{(2\sigma_1 - 1)\log q} < \varepsilon_1,
\]
where the last inequality follows from our explicit choice of $m > \ell_0$. For the third sum on the right-hand side of~\eqref{eq_m_to_hg}, each term has $n \deg P > h_g$, so for all $s \in K$ we have
\[
\Bigg| \sum_{\substack{P\ \text{prime} \\ \deg P \le m}} \sum_{n > \lfloor h_g / \deg P \rfloor} \frac{a_P^n}{n |P|^{ns}} \Bigg|
\le \sum_{\substack{P\ \text{prime} \\ \deg P \le m}} \frac{q^{-h_g \sigma_1}}{1 - q^{-\sigma_1}} < \varepsilon_1
\]
for all sufficiently large $g$, because $m$ is fixed and $h_g \to \infty$ as $g \to \infty$. Combining these three estimates via the triangle inequality, we conclude that for all sufficiently large $g$ and every $D \in \mathcal{A}_{2g+1}^m$,
\[
\max_{s \in K} \Bigg| \sum_{\substack{P\ \text{prime} \\ \deg P \le h_g}} \sum_{n=1}^{\lfloor h_g / \deg P \rfloor} \frac{\chi_D(P^n)}{n |P|^{ns}} - \sum_{\substack{P\ \text{prime} \\ \deg P \le m}} \sum_{n=1}^{\infty} \frac{a_P^n}{n |P|^{ns}} \Bigg| < 3\varepsilon_1.
\]
Together with~\eqref{eq_fund_lemma_approx}, this bound gives us
\begin{equation}\label{eq_f_diff_hg}
\max_{s \in K} \Bigg| f(s) - \sum_{\substack{P\ \text{prime} \\ \deg P \le h_g}} \sum_{n=1}^{\lfloor h_g / \deg P \rfloor} \frac{\chi_D(P^n)}{n |P|^{ns}} \Bigg| < 4\varepsilon_1.
\end{equation}

Having approximated $f(s)$ by a sum over primes $P$ with degree $\le h_g$, we now compare $F(s) = \exp(f(s))$ with the truncated product $P_k(s, \chi_D)$ defined in Lemma~\ref{lem_trunc_approx} with $k = h_g$, namely
\[
P_{h_g}(s, \chi_D) = \exp\Bigg( \sum_{\substack{P\ \text{prime} \\ \deg P \le h_g}} \sum_{n=1}^{\lfloor h_g / \deg P \rfloor} \frac{\chi_D(P^n)}{n |P|^{ns}} \Bigg).
\]
Recall that $\exp(f(s)) = F(s)$ and $\displaystyle M = \max_{s \in K} |F(s)|$. For any $s \in K$ and $D \in \mathcal{A}_{2g+1}^m$, factoring out $F(s) = \exp(f(s))$ gives
\begin{align}\label{eq_exp_factor}
|P_{h_g}(s, \chi_D) - F(s)|
&= |F(s)| \cdot \Bigg| \exp\Bigg( \sum_{\substack{P\ \text{prime} \\ \deg P \le h_g}} \sum_{n=1}^{\lfloor h_g / \deg P \rfloor} \frac{\chi_D(P^n)}{n |P|^{ns}} - f(s) \Bigg) - 1\Bigg| \notag \\
&< 8M\varepsilon_1,
\end{align}
where the last inequality follows from~\eqref{eq_choice_of_epsilon1}, \eqref{eq_f_diff_hg}, and the elementary bound $|\exp(w) - 1| \le 2|w|$ valid for $|w| \le \frac12$. In particular, by the triangle inequality and~\eqref{eq_choice_of_epsilon1}, we obtain the uniform upper bound
\begin{equation}\label{eq_Phg_upper}
\max_{s \in K} |P_{h_g}(s, \chi_D)| \le \max_{s \in K} |F(s)| + \max_{s \in K} |P_{h_g}(s, \chi_D) - F(s)| \le M + 8M\varepsilon_1 \le 2M.
\end{equation}

Next, we apply Lemma~\ref{lem_trunc_approx} with $k = h_g$ to pass from the truncated product $P_{h_g}(s, \chi_D)$ to $L(s, \chi_D)$ itself. By our choice of $h_g = \left\lfloor \frac{2g}{(2\sigma_1 - 1)\log_q g} \right\rfloor$, we see that
\[
\lim_{g \to \infty} \frac{g}{h_g} q^{\left(\frac12 - \sigma_1\right)h_g} = 0.
\]
In particular, $\frac{g}{h_g} q^{\left(\frac12 - \sigma_1\right)h_g} \le 1$ for all sufficiently large $g$, so~\eqref{eq_hybrid_exp} of Lemma~\ref{lem_trunc_approx} gives, uniformly for all $s \in K$,
\[
|L(s, \chi_D) - P_{h_g}(s, \chi_D)| = |P_{h_g}(s, \chi_D)| \cdot O\left( \frac{g}{h_g} q^{\left(\frac12 - \sigma_1\right)h_g} \right),
\]
where the implied constant depends at most on $q$ and $\sigma_1$. Combining this limit with the uniform bound $\max_{s \in K} |P_{h_g}(s, \chi_D)| \le 2M$ from~\eqref{eq_Phg_upper}, we find that for all sufficiently large $g$ and every $D \in \mathcal{A}_{2g+1}^m$,
\begin{equation}\label{eq_L_Phg_bound}
\max_{s \in K} |L(s, \chi_D) - P_{h_g}(s, \chi_D)| < \varepsilon_1.
\end{equation}

Finally, we combine the relevant approximations above. For every $D \in \mathcal{A}_{2g+1}^m$ and all sufficiently large $g$, using the triangle inequality together with~\eqref{eq_exp_factor} and~\eqref{eq_L_Phg_bound} gives
\[
\max_{s \in K} |L(s, \chi_D) - F(s)| < \varepsilon_1 + 8M\varepsilon_1 < \varepsilon,
\]
where the last inequality follows from our choice of $\varepsilon_1$ in~\eqref{eq_choice_of_epsilon1}. Because every polynomial $D \in \mathcal{A}_{2g+1}^m$ satisfies this inequality for all sufficiently large $g$, it follows from the density bound~\eqref{eq_positive_density} that
\begin{align*}
&\liminf_{g \to \infty} \frac{1}{\#\Hg} \# \left\{ D \in \Hg : \max_{s \in K} |L(s, \chi_D) - F(s)| < \varepsilon \right\} \\
&\ge \liminf_{g \to \infty} \frac{\# \mathcal{A}_{2g+1}^m}{\#\Hg} \ge \frac12 C_m > 0.
\end{align*}
This completes the proof.

\section*{Acknowledgements}

I would like to express my sincere gratitude to Professor Steven M.~Gonek for suggesting this problem and for providing invaluable guidance, encouragement, and support throughout the entire duration of this project. I am also deeply grateful to all my teachers and professors who have guided, mentored, and supported me at every stage of my academic journey and life.

\end{document}